\documentclass[preprint,10pt]{elsarticle}

\usepackage{amsmath,amsthm,amssymb,color}
\usepackage{graphicx}
\usepackage{comment}

\def\a{{\bf a}}
\def\b{{\bf b}}
\def\c{{\bf c}}
\def\d{{\bf d}}

\def\f{{\bf f}}
\def\g{{\bf g}}
\def\h{{\bf h}}

\def\p{{\bf p}}
\def\q{{\bf q}}
\def\r{{\bf r}}
\def\s{{\bf s}}
\def\bu{{\bf u}}
\def\bv{{\bf v}}

\def\w{{\bf w}}

\def\z{{\bf z}}

\def\rd{{\rm d}}
\def\re{{\rm e}}
\def\ri{{\rm i}}

\newcommand{\pc}{\overline{\p}}
\newcommand{\qc}{\overline{\q}}
\newcommand{\sbar}{\overline{\s}}
\newcommand{\buc}{\overline{\bu}}
\newcommand{\bvc}{\overline{\bv}}

\def\AMM{{\it Amer.\ Math.\ Monthly\ }}
\def\ACM{{\it Adv.\ Comp.\ Math.\ }}
\def\ACMTMS{{\it ACM Trans.\ Math.\ Software\ }}
\def\ACMTOG{{\it ACM Trans.\ Graphics\ }}
\def\CAD{{\it Comput.\ Aided Design }}
\def\CAEJ{{\it Comput.\ Aided Eng. J.\ }}
\def\CAGD{{\it Comput.\ Aided Geom.\ Design }}
\def\CG{{\it Computers \& Graphics }}
\def\CVGIP{{\it Comput.\ Vision, Graphics, Image\ Proc.\ }}
\def\IBMJRD{{\it IBM J.\ Res.\ Develop.\ }}

\def\JMAA{{\it J.\ Math.\ Anal.\ Appl.\ }}
\def\MC{{\it Math.\ Comp.\ }}
\def\MMAS{{\it Math.\ Methods\ Appl.\ Sci.\ }}
\def\NA{{\it Numer.\ Algor.\ }}
\def\PAMS{{\it Proc.\ Amer.\ Math.\ Soc.\ }}
\def\SIAMJNA{{\it SIAM J.\ Numer.\ Anal.\ }}
\def\SIAMR{{\it SIAM Rev.\ }}

\newcommand{\be}{\begin{equation}}
\newcommand{\ee}{\end{equation}}
\newcommand{\ba}{\begin{eqnarray}}
\newcommand{\ea}{\end{eqnarray}}
\newcommand{\bi}{\begin{itemize}}
\newcommand{\ei}{\end{itemize}}

\newtheorem{rmk}{Remark}
\newtheorem{lma}{Lemma}

\newtheorem{exm}{Example}

\newcommand{\RR}{\mathbb{R}}
\newcommand{\CC}{\mathbb{C}}
\newcommand{\VCR}{\boldsymbol{C}_{\!R}}
\newcommand{\VCI}{\boldsymbol{C}_{\!I}}

\begin{document}
\begin{frontmatter}
    
\title{
Application of a metric for complex 
polynomials to bounded modification of 
planar Pythagorean--hodograph curves
}

\author[address1]{Rida T. Farouki}
\ead{farouki@ucdavis.edu}

\author[address2,address3]{Marjeta Knez\corref{corauth}}
\ead{marjetka.knez@fmf.uni-lj.si}
\cortext[corauth]{Corresponding author}

\author[address4,address5]{Vito Vitrih}
\ead{vito.vitrih@upr.si}

\author[address2,address3]{Emil \v{Z}agar}
\ead{emil.zagar@fmf.uni-lj.si}

\address[address1]{Mechanical \& Aerospace Engineering, 
University of California, Davis, CA 95616, USA}
\address[address2]{Faculty of Mathematics and Physics, University of Ljubljana, Jadranska 19, Ljubljana, Slovenia}
\address[address3]{Institute of Mathematics, Physics and Mechanics, Jadranska 19, Ljubljana, Slovenia}
\address[address4]{Faculty of Mathematics, 
Natural Sciences and Information Technologies, \\
University of Primorska, Glagolja\v{s}ka 8, Koper, Slovenia}
\address[address5]{Andrej Maru\v{s}i\v{c} Institute,
University of Primorska, Muzejski trg 2, Koper, Slovenia}

\begin{abstract}
\noindent
By interpreting planar polynomial curves as complex--valued functions of a 
real parameter, an inner product, norm, metric function, and the notion of 
orthogonality may be defined for such curves. This approach is applied to the 
complex pre--image polynomials that generate planar Pythagorean--hodograph 
(PH) curves, to facilitate the implementation of bounded modifications of them 
that preserve their PH nature. The problems of bounded modifications under the 
constraint of fixed curve end points and end tangent directions, and of increasing the 
arc length of a PH curve by a prescribed amount, are also addressed. 
\end{abstract}

\end{frontmatter}

\section{Introduction}

In the complex model \cite{farouki94} for planar PH curves, points $(x,y)$ 
in the Euclidean plane are identified with the complex values $x+\ri\,y$. 
A planar PH curve $\r(t)$ for $t\in[\,0,1]$ may be generated from a complex
\emph{pre--image polynomial} $\w(t)$ by integrating the
hodograph expression $\r'(t)=\w^2(t)$. This guarantees that the components
of $\r'(t)=x'(t)+\ri\,y'(t)$ satisfy \cite{farouki90} the polynomial
Pythagorean condition
\[
x'^2(t)+y'^2(t) \,=\, \sigma^2(t) \,,
\]
where $\sigma(t)=|\r'(t)|=|\w(t)|^2$ is the \emph{parametric speed} of $\r(t)$
--- the derivative $\rd s/\rd t$ of the curve arc length $s$ with respect to
the parameter $t$.

Planar PH curves admit an exact computation of quantities such as arc 
lengths and offset curves \cite{farouki08}, that necessitate use of numerical 
approximation for ``ordinary'' polynomial curves. However, their non--linear 
nature entails more sophisticated construction algorithms, and renders \emph
{a posteriori} shape modification a difficult task. To address this latter
problem, it is important to first formulate a measure of ``how close'' two
PH curves are, i.e., to specify a \emph{metric} for the space of all planar
PH curves.

The complex representation of planar PH curves offers a solution to this
task, in terms of the standard concepts of inner products and norms from
functional analysis \cite{kreyszig}. By introducing a bound on the distance
between an original and modified pre--image polynomial, it is possible to
characterize the set of changes to its coefficients that define the shape
modifications to a planar PH curve that do not compromise its PH nature. 

The focus of the methodology presented herein is on the planar PH curves, 
although the approach may be adaptable to the spatial PH curves \cite
{choi02,farouki02} and the numerous other formulations of PH curves with 
distinctive properties that have been proposed \cite
{aithaddou17,albrecht17,kim17,kim19,kosinka06,kosinka10,moon99,pottmann95a,pottmann95b,romani19,romani14}.

The remainder of this paper is organized as follows. Section~\ref{sec:metric}
introduces the basic concepts of an \emph{inner product}, \emph{norm}, and 
\emph{metric} for the space of all polynomials in a real variable $t\in
[\,0,1\,]$ with complex coefficients. Section~\ref{sec:ocurves} then shows 
that this metric allows an \emph{angle} between such polynomials to be 
defined, and gives examples of \emph{orthogonal} plane curves specified as
complex polynomial functions of a real parameter. Section~\ref{sec:phcurves}
discusses the application of these concepts to planar PH curves and it is
observed that to maintain the PH nature, modifications should be made to 
the pre-image polynomial instead of directly to the curve. Modifications 
satisfying a prescribed bound on the distance between the original and
modified pre--image polynomials are discussed in Section~\ref{sec:pre-image},
and modifications that preserve the end tangents or end points of PH curves
are also formulated. Section~\ref{sec:arclen} shows how complex polynomials 
orthogonal to a specified pre--image polynomial may be used to modify PH 
curve arc lengths. Finally, Section~\ref{sec:close} summarizes the 
contributions of this study and suggests further possible avenues of 
investigation.

\section{Metric space of complex polynomials}
\label{sec:metric}

We begin by reviewing some elementary concepts from functional analysis
--- inner products, norms, and metrics (see \cite{kreyszig} for a thorough 
treatment). 

Let $\bu(t),\bv(t)\in \mathbb{C}[t]$ be complex polynomials in the real
variable $t\in[\,0,1\,]$. Their complex--valued \emph{inner product} is 
defined by
\be
\langle\bu,\bv\rangle \,=\, \int_0^1 \bu(t)\,\bvc(t) \; \rd t \,.
\nonumber
\ee
For any complex polynomial $\w(t)\in\mathbb{C}[t]$, this inner product induces
a non--negative norm specified by
\be
\label{norm}
\|\w\| \,=\, \sqrt{\langle\w,\w\rangle} \,.
\ee

A metric, or \emph{distance function}, for the complex polynomials $\bu(t)$
and $\bv(t)$ may be defined in terms of the norm (\ref{norm}) as
\be
\label{metric}
\mbox{distance}(\bu,\bv) \,=\, \|\bu-\bv\| \,.
\ee
Since
\ba
\|\bu-\bv\|^2 \!\! &=& \!\!
\int_0^1 (\bu(t)-\bv(t))(\buc(t)-\bvc(t)) \, \rd t 
\nonumber \\
\!\! &=& \!\!
\int_0^1 |\bu(t)|^2+|\bv(t)|^2-2\,\mbox{Re}(\bu(t)\bvc(t)) \, \rd t 
\nonumber \\
\!\! &=& \!\!
\|\bu\|^2+\|\bv\|^2-2\,\mbox{Re}(\langle\bu,\bv\rangle) \,,
\nonumber
\ea
we have
\[
\mbox{distance}(\bu,\bv) \,=\,
\sqrt{\|\bu\|^2+\|\bv\|^2-2\,\mbox{Re}(\langle\bu,\bv\rangle)} \,.
\]
Note that $\mbox{distance}(\bu,\bv)=0$ if and only if $\bu(t)\equiv\bv(t)$.

The metric (\ref{metric}) can be used to define the distance between 
planar curves, $\r(t)$ and $\s(t)$, regarded as complex functions of the 
real parameter $t\in[\,0,1\,]$ --- namely,
\[
\mbox{distance}(\r,\s) \,=\,
\sqrt{\|\r\|^2+\|\s\|^2-2\,\mbox{Re}(\langle\r,\s\rangle)} \,.
\]
Note that, when $\r(t)$ and $\s(t)$ are orthogonal, i.e., $\mbox{Re}
(\langle\r,\s\rangle)=0$, the distance is simply $\sqrt{\|\r\|^2+\|\s\|^2}$. 
The following elementary cases are noteworthy.
\begin{itemize}
\item[1.]
If $\s(t)$ is a translate of $\r(t)$ by  the complex value $\d$, $\mbox
{distance}(\r,\s)=|\d|$.
\item[2.]
If $\s(t)$ is a rotation of $\r(t)$ by angle $\theta$, $\mbox{distance}
(\r,\s)=2\,(1-\cos\theta)\,\|\r\|$.
\item[3.]
If $\s(t)$ is a scaling of $\r(t)$ by the factor $c$, $\mbox{distance}
(\r,\s)=|1-c\,|\,\|\r\|$.
\end{itemize}

In some contexts it may be desirable for $\mbox{distance}(\r,\s)$ to reflect 
only the differences of \emph{shape}, and discount considerations of position, 
orientation, and scaling. If $\r(0)=\s(0)=0$, this can be achieved through 
a rotation/scaling transformation that makes the vectors $\r(1)-\r(0)$ and 
$\s(1)-\s(0)$ coincident.

The preceding ideas were briefly mentioned in the problem of constructing
spatial $C^2$ closed loops with prescribed arc lengths using PH curves \cite
{farouki21} --- the solutions can be characterized in terms of two complex 
polynomials $\bu(t),\bv(t)$ satisfying $\|\bu\|=\|\bv\|=1/\sqrt{2}$, $\langle
\bu,\bv\rangle=0$, and thus $\mbox{distance}(\bu,\bv)=1$.

\section{Orthogonal planar curves}
\label{sec:ocurves}

Since $|\,\mbox{Re}(\langle\bu,\bv\rangle)\,| \le \|\bu\|\,\|\bv\|$ from the 
Cauchy--Schwartz inequality, an angle $\theta\in[\,0,\pi\,]$ between $\bu$ 
and $\bv$ may be defined by
\be
\cos\theta \,=\, \frac{\mbox{Re}(\langle\bu,\bv\rangle)}{\|\bu\|\,\|\bv\|} \,,\nonumber
\ee
and we thereby obtain the cosine rule
\[
\mbox{distance}^2(\bu,\bv) \,=\, \|\bu\|^2+\|\bv\|^2
- 2\,\|\bu\|\,\|\bv\|\cos\theta \,.
\]
If $\cos\theta=0$ --- i.e., $\mbox{Re}(\langle\bu,\bv\rangle)=0$ --- we say 
that $\bu$ and $\bv$ are \emph{orthogonal}, and write $\bu\perp\bv$.
For two orthogonal complex polynomials, the distance becomes simply $\sqrt
{\|\bu\|^2+\|\bv\|^2}$.

Although the focus herein is on PH curves, the above principles apply to 
\emph{any} planar curves represented as complex-valued polynomial functions 
of a real variable, an approach to the study of planar curves promoted by
Zwikker \cite{zwikker}. If $\r(t)$ and $\s(t)$ are B\'ezier curves of degree 
$m$ and $n$, with control points $\p_0,\ldots,\p_m$ and $\q_0,\ldots,\q_n$, 
the product $\r(t)\,\overline{\s}(t)$ can be expressed \cite{farouki12} as
\[
\r(t)\,\overline{\s}(t) \,=\, \sum_{k=0}^{m+n} \z_k 
\binom{m+n}{k}(1-t)^{m+n-k}t^k \,,
\]
with
\[
\z_k \,=\, \sum_{j=\max(0,k-n)}^{\min(m,k)}
\frac{\displaystyle\binom{m}{j}\binom{n}{k-j}}
{\displaystyle\binom{m+n}{k}}
\, \p_j\,\qc_{k-j} \,, \quad k=0,\ldots,m+n \,.
\]
Since the definite integral of every Bernstein basis function of degree
$m+n$ over $[\,0,1\,]$ is simply $1/(m+n+1)$, the inner product of $\r(t)$
and $\s(t)$ is
\[
\langle\r,\s\rangle \,=\, \frac{\z_0+\cdots+\z_{m+n}}{m+n+1} \,.
\]
Thus, for given control points $\p_0,\ldots,\p_m$ of $\r(t)$, the orthogonality
condition $\mbox{Re}(\langle\r,\s\rangle)=0$ amounts to a single linear 
constraint on the real and imaginary parts of the control points $\q_0,
\ldots,\q_n$ of $\s(t)$, so the dimension of the subspace of degree $n$ 
curves $\s(t)$ that are orthogonal to $\r(t)$ is $2\,n+1$. To explore this
subspace in more detail, we set
\[
(r_x(t),r_y(t)) \,=\, (\mbox{Re}(\r(t)),\mbox{Im}(\r(t))) \,, \;\;
(s_x(t),s_y(t)) \,=\, (\mbox{Re}(\s(t)),\mbox{Im}(\s(t))) \,,
\] 
and define
\begin{align*}
d(t) \,:=\,
\mbox{Re}(\r(t)\,\sbar(t)) & \,=\, r_x(t)s_x(t)+r_y(t)s_y(t) \,, \\
c(t) \,:=\,
\mbox{Im}(\r(t)\,\sbar(t)) & \,=\, s_x(t)r_y(t)-s_y(t)r_x(t) \,.
\end{align*}
Regarding $\r(t),\s(t)$ as vector functions, $d(t)$ is their dot product 
and $c(t)$ is the component of the cross product orthogonal to the $(x,y)$ 
plane. Moreover, $\mbox{Re}(\langle\r,\s\rangle)$ and $\mbox{Im}(\langle\r,
\s\rangle)$ are the integrals of $d(t)$ and $c(t)$ over $t\in[\,0,1\,]$.

To construct orthogonal curves, it is convenient to employ an orthonormal
polynomial basis. We choose here the Legendre basis on $t\in[\,0,1\,]$ which 
may be expressed in terms of the Bernstein basis \cite{farouki00} as 
\[
L_k(t) \,=\, \sqrt{2k+1}\, 
\sum_{i=0}^k (-1)^{k+i} \binom{k}{i}\,b^k_i(t) \,,
\quad b^k_i(t) \,=\, \binom{k}{i}(1-t)^{k-i}t^i \,.
\]
These basis functions satisfy 
\be
\int_0^1 L_j(t)L_k(t) \; \rd t \,=\, \delta_{jk} \,,\nonumber
\ee
where $\delta_{jk}$ is the Kronecker delta, and the first few instances are
\ba
L_0(t) \!\! &=& \!\!  1 \,,
 \nonumber \\
L_1(t) \!\! &=& \!\! \sqrt{3}\,(2\,t-1) \,,
\nonumber \\
L_2(t) \!\! &=& \!\! \sqrt{5}\,(6\,t^2-6\,t+1) \,,
 \nonumber \\
L_3(t) \!\! &=& \!\! \sqrt{7}\,(20\,t^3-30\,t^2+12\,t-1) \,.
\nonumber
\ea

For any given curve $\r(t)$, $t\in[\,0,1\,]$ of degree $m$ we consider the
problem of constructing curves $\r_\perp(t)$, of the same degree, that are
orthogonal to $\r(t)$. Expressing $\r(t)$ and $\r_\perp(t)$ in the Legendre
basis as
\ba
\r(t) \!\! &=& \!\!
\sum_{k=0}^m a_{k,1} L_k(t) \,+\, \ri \sum_{k=0}^m a_{k,2} L_k(t) \,,
\nonumber \\
\r_\perp(t) \!\! &=& \!\!
\sum_{k=0}^m b_{k,1} L_k(t) \,+\, \ri \sum_{k=0}^m b_{k,2} L_k(t) \,,
\label{orth-curve-1}
\ea
by the orthonormality of the basis functions we have
\begin{align*}
\mbox{Re}(\langle\r,\r_\perp\rangle) & \,=\, 
\int_{0}^1 \sum_{k=0}^m a_{k,1} L_{k}(t) \sum_{\ell=0}^m b_{\ell,1} L_{\ell}(t) +
\sum_{k=0}^m a_{k,2} L_{k}(t) \sum_{\ell=0}^m b_{\ell,2} L_{\ell}(t) \, \rd t \\
& \,=\, \sum_{k=0}^m a_{k,1}b_{k,1} + a_{k,2}b_{k,2} \,.
\end{align*}
Thus, identifying the coefficients $b_{k,1}+\ri\,b_{k,2}$, $k=0,1,\dots,m$, 
of an orthogonal curve is equivalent to finding the set of $2m+1$ linearly 
independent vectors
$$
\boldsymbol{b} \,=\, (b_{0,1},b_{0,2},b_{1,1},b_{1,2},\dots,b_{m,1},b_{m,2})^T
\in \mathbb{R}^{2m+2}
$$ 
that are orthogonal to the vector 
$$
\boldsymbol{a} \,=\, (a_{0,1},a_{0,2},a_{1,1},a_{1,2},\dots,a_{m,1},a_{m,2})^T 
\in \mathbb{R}^{2m+2}
$$ 
with respect to the Euclidean inner (or dot) product in $\mathbb{R}^{2m+2}$. 
The basis of the orthogonal complement $\boldsymbol{a}^{\perp}$ follows from 
the extended QR decomposition $\boldsymbol{a}=QR$, where $Q$ is a $(2m+2)
\times(2m+2)$ orthogonal matrix, and $R=(\left\|\boldsymbol{a}\right\|_2,0,
\dots,0)^T$ where $\|\cdot\|_2$ is the standard Euclidean norm. The matrix 
$Q$ is the well--known Householder reflection \cite
{hamming}, computed as
\be
\label{Q}
Q \,=\, I - \frac{2}{\boldsymbol{g}^T\!\boldsymbol{g}}\,
\boldsymbol{g}\boldsymbol{g}^T \,, \quad 
\boldsymbol{g} \,=\, \boldsymbol{a} + \mbox{sign}(a_{0,1}) 
\left\|\boldsymbol{a}\right\|_2 \boldsymbol{e}_1 \,,
\ee
where $\boldsymbol{e}_1=(1,0,\dots,0)^T$. The second through last columns of 
$Q$ --- denoted by $\boldsymbol{q}_2,\dots,\boldsymbol{q}_{2m+2}$ --- are the 
orthogonal basis of $\boldsymbol{a}^{\perp}$. Thus, any vector 
$$
\boldsymbol{b} \,=\, \sum_{k=2}^{2m+2} \xi_{k-1}\boldsymbol{q}_k \,, \quad 
\xi_{1},\dots,\xi_{2m+1}\in\mathbb{R}
$$
identifies a curve $\r_\perp(t)$ of the form \eqref{orth-curve-1}, that is 
orthogonal to $\r(t)$. Moreover, the columns of $Q$ also define curves that 
are pairwise orthogonal.

\begin{exm}
\label{exm:orthogonal-1}
{\rm Consider the vector $\boldsymbol{a}=\left(\alpha_0,0,0,\alpha_1,\alpha_2,
0,0,\alpha_3\right)$ that defines the curve 
$$
\r(t) \,=\, \alpha_0L_0(t)+\alpha_2L_2(t)
+ \ri\,(\alpha_1L_1(t)+\alpha_3L_3(t)) \,,
$$ 
with B\'ezier control points 
$$
\p_0 \,=\, \alpha_0+\sqrt{5}\,\alpha_2 - 
\ri\,(\sqrt{3}\,\alpha_1 +\sqrt{7}\,\alpha_3) \,=\, \pc_3 \,,
$$
$$
\p_1 \,=\, \alpha_0-\sqrt{5}\,\alpha_2 - 
\ri\,\left(\frac{\alpha_1}{\sqrt{3}}- 3\sqrt{7}\,\alpha_3\right) 
\,=\  \pc_2 \,.
$$
In this case, we obtain 
\[
Q= \left[
\begin{array}{cccccccc}
-\frac{\alpha_0}{\alpha } & 0 & 0 & -\frac{\alpha_1}{\alpha } & 
-\frac{\alpha_2}{\alpha } & 0 & 0 & -\frac{\alpha_3}{\alpha } \\
0 & 1 & 0 & 0 & 0 & 0 & 0 & 0 \\
0 & 0 & 1 & 0 & 0 & 0 & 0 & 0 \\
-\frac{\alpha_1}{\alpha} & 0 & 0 & 1-\frac{\alpha_1^2}{\alpha^2+\alpha_0
\alpha} & -\frac{\alpha_1\alpha_2}{\alpha^2+\alpha_0\alpha} & 0 & 0 &
-\frac{\alpha_1\alpha_3}{\alpha^2+\alpha_0\alpha} \\
-\frac{\alpha_2}{\alpha} & 0 & 0 & -\frac{\alpha_1\alpha_2}
{\alpha^2+\alpha_0\alpha} & 1-\frac{\alpha_2^2}{\alpha^2+\alpha_0\alpha} 
& 0 & 0 & -\frac{\alpha_2\alpha_3}{\alpha^2+\alpha_0\alpha} \\
0 & 0 & 0 & 0 & 0 & 1 & 0 & 0 \\
0 & 0 & 0 & 0 & 0 & 0 & 1 & 0 \\
-\frac{\alpha_3}{\alpha} & 0 & 0 & -\frac{\alpha_1\alpha_3}{\alpha^2
+\alpha_0\alpha} & -\frac{\alpha_2\alpha_3}{\alpha^2+\alpha_0\alpha} & 
0 & 0 & 1-\frac{\alpha_3^2}{\alpha ^2+\alpha_0\alpha} \\
\end{array}
\right],
\]
where $\alpha=\|\boldsymbol{a}\|_2$. From columns $2,3,6,7$, we see that any 
curve of the form 
\begin{equation}
\label{orthogonalCurve2}
\r_\perp(t) \,=\,
\beta_1 L_1(t)+\beta_3 L_3(t)+\ri\,(\beta_0 L_0(t)+ \beta_2 L_2(t))    
\end{equation}
is orthogonal to $\r(t)$. The B\'ezier control points of such a curve are
$$
\q_0 \,=\, -\sqrt{3}\,\beta_1-\sqrt{7}\,\beta_3
+ \ri\,(\beta_0+\sqrt{5}\,\beta_2) \,=\, -\,\qc_3 \,,
$$
$$
\q_1 \,=\, -\frac{\beta_1}{\sqrt{3}}+3\sqrt{7}\,\beta_3 + 
\ri(\beta_0-\sqrt{5}\,\beta_2) \,=\, -\,\qc_2 \,.
$$
Columns $4,5,8$ of $Q$ define three additional linearly--independent 
orthogonal curves, that are symmetric about the real axis. 

As an illustrative example, consider the case 
$$
\alpha_0=1, \,\, \alpha_1=-2\sqrt{3}, \,\, 
\alpha_2=\sqrt{5}, \,\, \alpha_3=\sqrt{7}.
$$
Figure~\ref{fig:example-orthogonal1} shows the curve $\r(t)$ (black) and 
four curves $\r_\perp(t)$ orthogonal to it, corresponding to the columns 4 
(red), 5 (green) and 8 (blue) of $Q$, multiplied by the norm of $\r(t)$,
while the purple curve is defined by \eqref{orthogonalCurve2} with 
$$
\beta_0=-1, \,\, \beta_1=3, \,\, \beta_1=-2, \beta_3=-4. 
$$
Figure~\ref{fig:example-orthogonal1} also shows the graphs of $\mbox{Re}
(\r(t)\,\r_\perp(t))$ for these four curves, which exhibit equal areas above 
and below the $t$--axis.
\begin{figure}[htb]
\begin{minipage}{0.49\textwidth}
\centering
\includegraphics[width=1\textwidth]{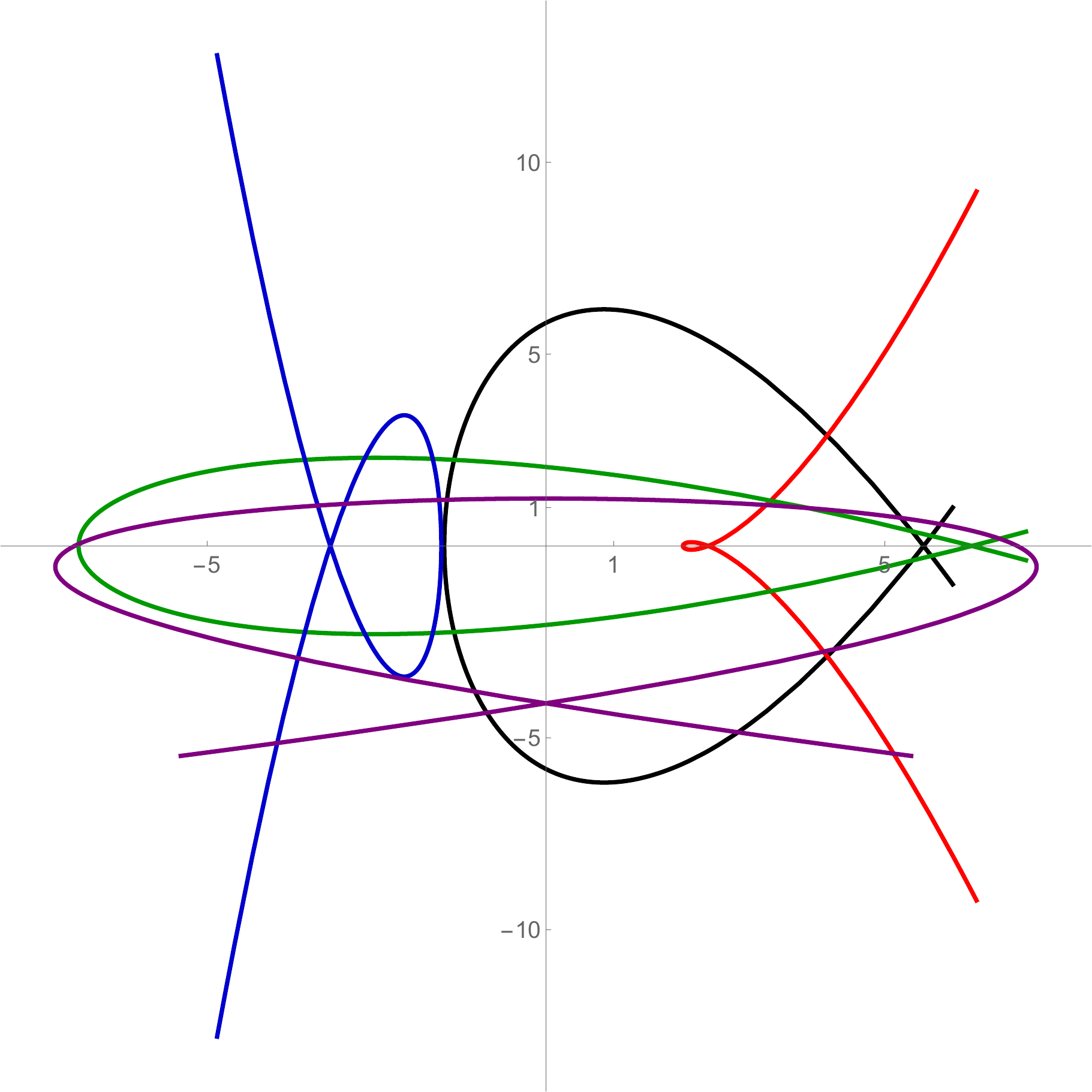}
\end{minipage}
\hskip2em
\begin{minipage}{0.4\textwidth}
\centering
\includegraphics[width=1\textwidth]{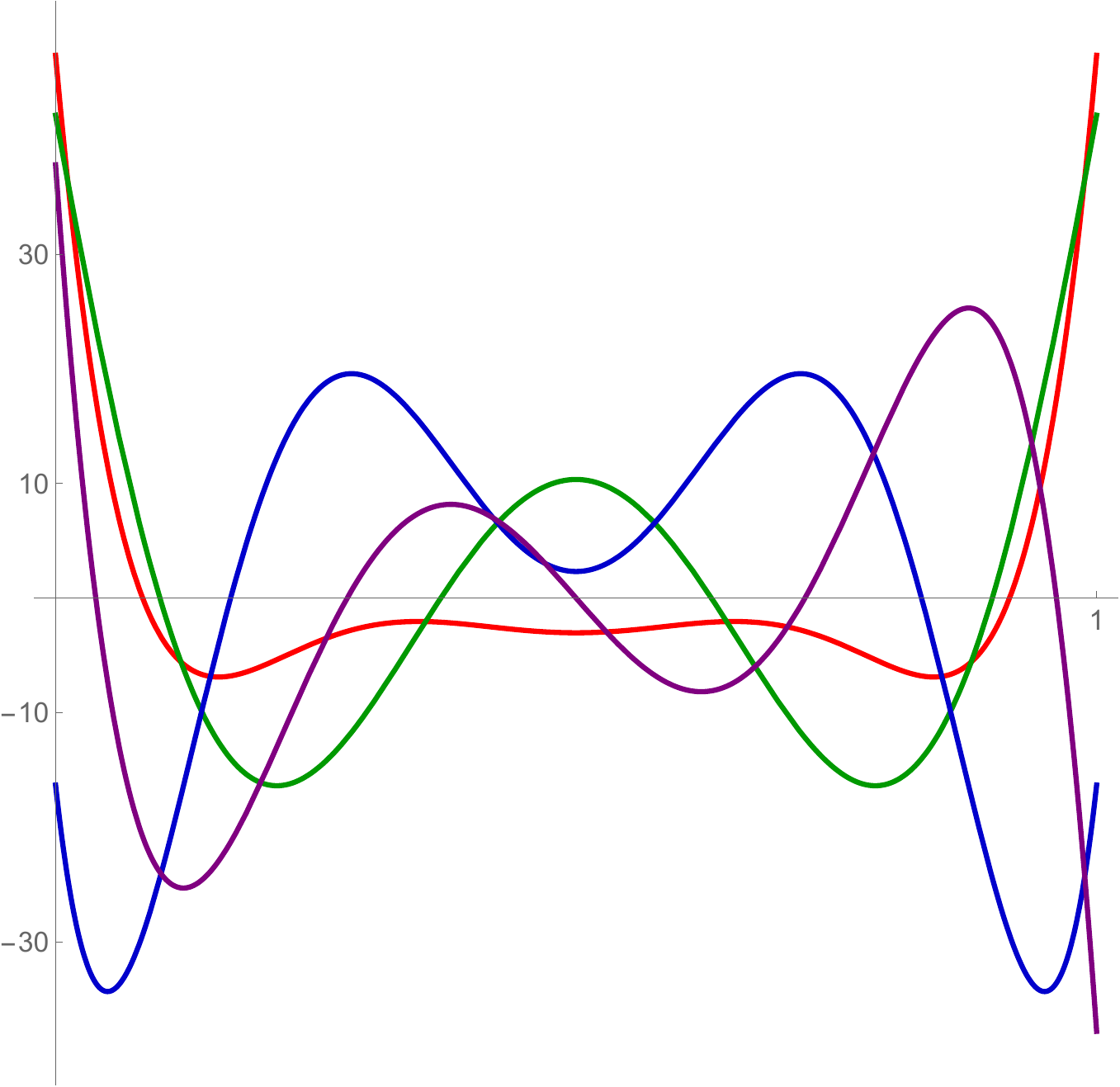}
\end{minipage}
\caption{Left: the curve $\r(t)$ in Example~\ref{exm:orthogonal-1} is
indicated in black, and four curves $\r_\perp(t)$ orthogonal to it are shown 
in red, green, blue, and purple. Right: the graphs of $\mbox{Re}(\r(t)\,
\r_\perp(t))$ for these four curves. }
\label{fig:example-orthogonal1}
\end{figure}
}
\end{exm}

\begin{exm}
\label{exm:orthogonal-2}
{\rm The cubic curve
$$
\r(t) \,=\, 7\,b^3_1(t) + \frac{16}{3}\,b^3_2(t) +  
\ri \left(\frac{20}{3}\,b^3_1(t) - \frac{11}{3}\,b^3_2(t) + 
\frac{19}{3}\,b^3_3(t)\right) 
$$
is a PH curve, since $\r'(t)=[\,5\,b^1_0(t)-3\,b^1_1(t)+\ri\left(2\,b^1_0(t)
-5\,b^1_1(t)\right)\,]^2$. From its Legendre coefficients, we obtain the vector 
$$
\boldsymbol{a} \,=\, 
\left(\frac{37}{12},\frac{7}{3},-\frac{\sqrt{3}}{12},\frac{13\sqrt{3}}{30},  
-\frac{37 \sqrt{5}}{60},\frac{\sqrt{5}}{6},\frac{\sqrt{7}}{28},
\frac{4\sqrt{7}}{15}\right) \,,
$$
and columns $2$--$8$ of its QR decomposition define $7$ orthogonal curves 
$\r_{\perp,k}(t)$, $k=1,\dots,7$. We compute their linear combination
$$
\r_{\perp}(t) \,=\, \sum_{k=1}^{7} \xi_k\,\r_{\perp,k}(t) \,,
$$ 
so that $\r_\perp(t)$ is a PH curve. To equate the number of equations and 
unknowns $\xi_1,\ldots,\xi_7$ we also require $\r_{\perp}(0)=0$, and $\r_{\perp}
(t)$ to have a prescribed parametric speed, $\sigma(t)=20-40\,t+38\,t^2$. 
The resulting non--linear system has six different solutions, illustrated in 
Figure~\ref{fig:example-orthogonal2}. Note that the prescribed $\sigma(t)$ 
implies that curves $\r(t)$ and $\r_{\perp}(t)$ all have the same arc length, 
namely $38/3$.

\begin{figure}[htbp]
\centering
\includegraphics[width=0.7\textwidth]{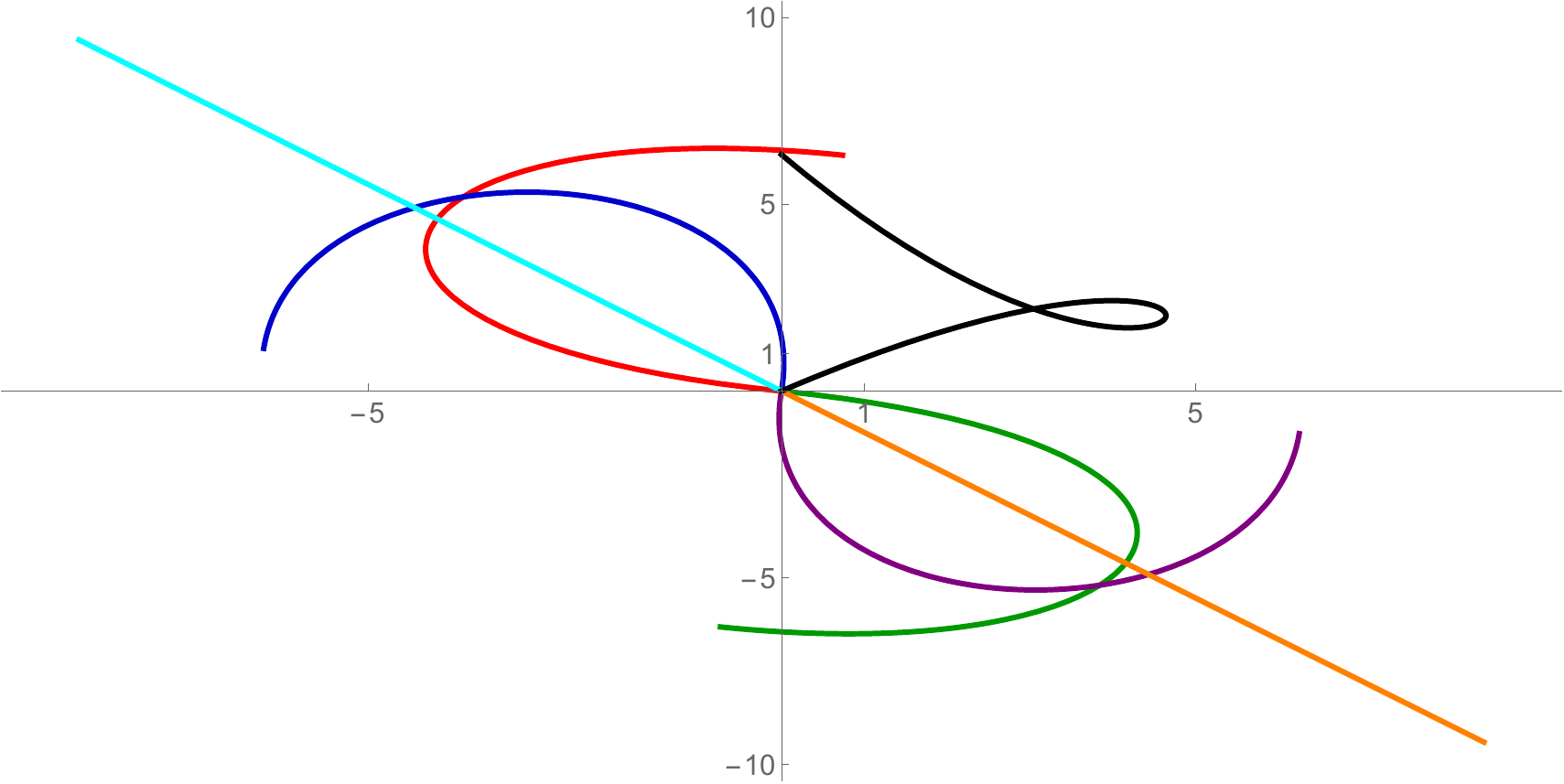}
\caption{The cubic PH curve $\r(t)$ (black) in Example~\ref{exm:orthogonal-2}, 
together with the six PH curves $\r_{\perp}(t)$ orthogonal to it (shown in 
different colors) that possess the same start point $(0,0)$ and have the prescribed parametric speed $\sigma(t)$.}
\label{fig:example-orthogonal2}
\end{figure}
}    
\end{exm}

\section{Planar Pythagorean-hodograph curves}
\label{sec:phcurves}

Planar PH curves $\r(t)$ are generated from complex pre--image polynomials 
$\w(t)$ by integrating the derivative or \emph{hodograph} $\h(t):=\r'(t)=
\w^2(t)$. If $\w(t)$ is of degree $m$, specified in Bernstein form as
\be
\label{w}
\w(t) \,=\, \sum_{k=0}^m \w_kb^m_k(t) \,,
\ee
the hodograph may be written as
\[
\h(t) \,=\, \sum_{k=0}^{2m} \h_kb^{2m}_k(t) \,,
\]
with coefficients determined \cite{farouki12} by
\be
\label{h_k}
\h_k \,=\, \sum_{j=\max(0,k-m)}^{\min(m,k)}
\frac{ \displaystyle\binom{m}{j}\binom{m}{k-j}}
{ \displaystyle\binom{2m}{k} }\,\w_j\w_{k-j} \,,
\quad 0 \le k \le 2m \,.
\ee
The B\'ezier control points of the PH curve of degree $n=2m+1$ constructed
by integrating $\r'(t)$ are then given by
\be
\label{ph-n}
\p_{k+1} \,=\, \p_k+\frac{\h_k}{2m+1} \,, \quad k=0,\ldots,n-1 \,,
\ee
where we henceforth assume $\p_0=0$ as the integration constant.
 
We focus mainly on the planar PH quintics, generated from a quadratic 
pre--image polynomial \eqref{w},  
which are widely considered to be the lowest--order PH curves appropriate to 
free--form design applications. The control points of the B\'ezier form 
\[
\r(t) \,=\, \sum_{k=0}^5 \p_kb^5_k(t) \,,
\]
are
\begin{equation*} 
\begin{split}
\p_1 & = \p_0 + \frac{1}{5}\,\w_0^2 \,, \quad \p_2 = \p_1 + \frac{1}{5}\,\w_0\w_1 \,, \\
\p_3  & = \p_2 + \frac{1}{5}\,\frac{2\,\w_1^2+\w_0\w_2}{3} \,, \\
\p_4 & =  \p_3 + \frac{1}{5}\,\w_1\w_2 \,, \quad \p_5 = \p_4 + \frac{1}{5}\,\w_2^2 \,.
\end{split}
\end{equation*}

Planar PH quintics are typically constructed as solutions to a first--order 
Hermite interpolation problem \cite{farouki95} for specified end points 
$\r(0),\r(1)$ and end derivatives $\r'(0)$, $\r'(1)$. It is not feasible to 
modify their shape \emph{a posteriori} by manipulating the control points, 
since this will ordinarily compromise their PH nature. Modifications that
preserve the PH property of a curve should be made to its pre--image
polynomial, rather than directly to the PH curve. The metric will therefore 
be primarily used to measure the distance between an original and modified 
pre--image polynomial of the PH curve it generates.

\section{Perturbation of pre-image polynomials}
\label{sec:pre-image}

A key application of the metric for complex polynomials is to provide a
means to make ``modest'' shape changes to PH curves that preserve the PH 
property. To achieve this, modifications must be made to the pre--image
polynomial. We consider perturbations $\delta\w(t)$ to a given pre--image 
polynomial $\w(t)$ that, for a prescribed bound $\Delta$, satisfy
\be
\label{ineq}
\mbox{distance}(\w,\w+\delta\w) \,=\, \|\delta\w\| \,\le\, \Delta \,.
\ee
The perturbed pre--image polynomial determines a perturbed PH curve, with
control points $\p_k$ displaced to $\p_k+\delta\p_k$ for $k=1,\ldots,n$. 
The perturbations $\delta\p_k$ may be obtained by replacing $\w_k$ by 
$\w_k+\delta\w_k$ in (\ref{h_k}) and (\ref{ph-n}), and they determine a 
perturbation $\delta\r(t)$, whose norm provides a measure of the difference 
between the modified and original curves.

Henceforth, we use the Legendre and Bernstein forms of both $\w(t)$ and
$\delta\w(t)$ interchangeably. Whereas the former offers more concise
formulations, the latter is the standard representation scheme in computer
aided geometric design and offers simpler implementation of certain
constraints, such as the preservation of initial/final tangent directions.

The following lemma \cite{farouki00} describes the transformation between
these two representations, which is known to be quite numerically stable.

\begin{lma} \label{lemmaCD}
For a polynomial $p(t)$ of degree $n$ expressed in the Legendre and Bernstein 
bases on $[\,0,1\,]$ as
$$
p(t) \,=\, \sum_{k=0}^n c_kL_k(t) \,=\, \sum_{j=0}^n d_j\,b_j^n(t) \,,
$$ 
the coefficients $\boldsymbol{C}=(c_0,\ldots,c_n)^T$ and $\boldsymbol{D}=
(d_0,\ldots,d_n)^T$ are related according to $\boldsymbol{D}=M_n\boldsymbol{C}$, 
where $M_n$ is the $(n+1)\times(n+1)$ matrix with elements
$$
M_{n,jk} \,=\, \frac{\sqrt{2k+1}}{\displaystyle\binom{n}{j}} 
\sum_{i=\max(0,j+k-n)}^{\min{j,k}} (-1)^{k+i} 
\binom{k}{i} \binom{k}{i} \binom{n-k}{j-i} \,,
\quad 0 \le j,k \le n \,,
$$
whose inverse $M_n^{-1}$ has elements
$$
M^{-1}_{n,jk} \,=\, \frac{\sqrt{2j+1}}{n+j+1} \binom{n}{k} 
\sum_{i=0}^j (-1)^{i+j} 
\frac{\displaystyle\binom{j}{i}\binom{j}{i}}{\displaystyle\binom{n+j}{k+i}} 
\,, \quad 0 \le j,k \le n \,.
$$
Note that the columns of the matrix $M_n$ are orthogonal.
\end{lma}

\begin{exm}
\label{exm:matrixM}
{\rm For the linear, quadratic, and cubic pre--image polynomials of cubic, 
quintic, and septic PH curves, the matrices $M_n$ and their inverses are
$$
M_1=\left[
\begin{array}{cc}
 1 & -\sqrt{3} \\[1mm]
 1 & \sqrt{3} 
\end{array}
\right], \quad 
M_1^{-1} = \left[
\begin{array}{cc}
 \frac{1}{2} & \frac{1}{2} \\[1mm]
 -\frac{\sqrt{3}}{6} & \frac{\sqrt{3}}{6} 
\end{array}
\right], \quad
$$
$$
M_2 = \left[
\begin{array}{ccc}
 1 & -\sqrt{3} & \sqrt{5} \\[1mm]
 1 & 0 & -2 \sqrt{5} \\[1mm]
 1 & \sqrt{3} & \sqrt{5} 
\end{array}
\right], \quad
M_2^{-1} = \left[
\begin{array}{ccc}
 \frac{1}{3} & \frac{1}{3} & \frac{1}{3} \\[1mm]
 -\frac{\sqrt{3}}{6} & 0 & \frac{\sqrt{3}}{6} \\[1mm]
 \frac{\sqrt{5}}{30} & -\frac{\sqrt{5}}{15} & \frac{\sqrt{5}}{30} 
\end{array}
\right], \quad
$$
$$
\quad M_3 = \left[
\begin{array}{cccc}
 1 & -\sqrt{3} & \sqrt{5} & -\sqrt{7} \\[1mm]
 1 & -\frac{\sqrt{3}}{3} & -\sqrt{5} & 3 \sqrt{7} \\[1mm]
 1 & \frac{\sqrt{3}}{3} & -\sqrt{5} & -3 \sqrt{7} \\[1mm]
 1 & \sqrt{3} & \sqrt{5} & \sqrt{7} 
\end{array}
\right], \quad
M_3^{-1} = \left[
\begin{array}{cccc}
\frac{1}{4} & \frac{1}{4} & \frac{1}{4} & \frac{1}{4} \\[1mm]
-\frac{3\sqrt{3}}{20} & -\frac{\sqrt{3}}{20} 
& \frac{\sqrt{3}}{20} & \frac{3 \sqrt{3}}{20} \\[1mm]
 \frac{\sqrt{5}}{20} & -\frac{\sqrt{5}}{20} & 
-\frac{\sqrt{5}}{20} & \frac{\sqrt{5}}{20} \\[1mm]
 -\frac{\sqrt{7}}{140} & \frac{\sqrt{7}}{140} 
& -\frac{3 \sqrt{7}}{140} & \frac{\sqrt{7}}{140} 
\end{array}
\right].
$$
}
\end{exm}
 
In the Legendre form, the pre--image polynomial (\ref{w}) is expressed
in terms of coefficients $\c_0,\ldots,\c_m$ as
\be
\label{wL-m}
\w(t) \,=\, \sum_{k=0}^{m}\c_k L_k(t)\,,
\ee
and the Bernstein coefficients can be recovered from the Legendre coefficients 
through the relations
\begin{equation}
\label{LB-m}
\w_j \,=\, \sum_{k=0}^m M_{m,jk}\,\c_k \,, \quad j=0,\ldots,m \,,
\end{equation}
from which the B\'ezier control points \eqref{ph-n} may be determined. The 
Legendre and Bernstein forms of the perturbation polynomial are
$$
\delta\w(t) \,=\, \sum_{k=0}^{m} \delta\c_k L_k(t) 
\; = \;
\sum_{j=0}^{m} \delta\w_j b_j^m(t)\,.
$$
The relations (\ref{LB-m}) also hold for $\delta\w_j$ in terms of 
$\delta\c_k$, i.e.,
\be
\label{vecWandC}
\delta\boldsymbol{W} \,=\, {M}_m \, \delta \boldsymbol{C}, \quad
\delta\boldsymbol{W} := (\delta\w_0,\dots,\delta\w_m)^T, \quad
\delta\boldsymbol{C} := (\delta\c_0,\dots,\delta\c_m)^T.
\ee
In the Legendre form, the norm of the perturbation $\delta\w(t)$ is
\be
\label{norm-deltaw-1}
\|\delta\w\| \,=\, \sqrt{|\delta\c_0|^2+\cdots+|\delta\c_m|^2} 
\,=\, \|\delta \boldsymbol{C}\|_2 \,,
\ee
where $\|\cdot\|_2$ again denotes the standard Euclidean norm. Writing perturbation coefficients as
\begin{equation} 
\label{eq:equalPer-L-gen}
\delta\c_k \,=\, \rho_k \exp(\ri\,\varphi_k), \quad k=0,\ldots,m \,,
\end{equation}
the inequality (\ref{ineq}) becomes simply
\be
\|\delta\w\| \,=\, \sqrt{\sum_{k=0}^m \rho_k^2} \,\le\, \Delta \,.
\nonumber
\ee
For perturbations of equal magnitude $\rho : = \rho_0 =\dots = \rho_m$
it further simplifies to 
\be
\|\delta\w\| \,=\, \sqrt{m+1}\,\rho \,\le\, \Delta \,,
\nonumber
\ee
so choosing $\rho\leq\Delta/\sqrt{m+1}$ satisfies the condition \eqref{ineq} 
for any given $\Delta$.

In the Bernstein form $\|\delta\w\|$ may be expressed, using \eqref{vecWandC} 
and \eqref{norm-deltaw-1}, as
\[
\|\delta\w\| \,=\, \|M_m^{-1} \, \delta \boldsymbol{W}\|_2.
\]
Since 
$$
\|M_m^{-1}\,\delta\boldsymbol{W}\|_2 \,\leq\, 
\|M_m^{-1} \|_2 \, \|\delta\boldsymbol{W}\|_2 \,,
$$
where $\|M_m^{-1}\|_2$ is the matrix norm induced by the Euclidean vector norm, 
which is equal to the largest singular value $\sigma_{\max}\left(M_m^{-1}\right) 
=1/\sqrt{m+1}$, any choice of coefficients $\delta \boldsymbol{W}$ such that 
$\|\delta\boldsymbol{W}\|_2\leq\sqrt{m+1}\,\Delta$ ensures satisfaction of
\eqref{ineq} for any given $\Delta$. 
If we express 
\begin{equation} 
\label{eq:equalPer-B-gen}
\delta\w_k \,=\, r_k \exp(\ri\, \phi_k) \,, \quad  k=0,\ldots,m \,,
\end{equation}
then $\|\delta\boldsymbol{W}\|_2 = \sqrt{\sum_{k=0}^m r_k^2}$, and in the case of equal--magnitude 
perturbations, $r:=r_0=\dots=r_m$, 
the simple choice $r \leq \Delta$ implies that \eqref{ineq} holds true. However,
this is just a sufficient condition. The inequality \eqref{ineq} can be 
satisfied for larger values of $r$ by using the inequalities
\begin{align}
& \|\delta\w\| \,=\, \frac{\displaystyle\sqrt{\Phi_{01}+2}}{\sqrt{3}}\,r 
\leq \Delta \,,
\nonumber \\
& \|\delta\w\| \,=\, \frac{\displaystyle\sqrt{3\,\Phi_{01}+\Phi_{02}
+3\,\Phi_{12}+8}}{\sqrt{15}}\,r \leq \Delta \,,
\label{eq-cond-for-m2} \\
& \|\delta\w\| \,=\, \frac{\displaystyle\sqrt{10\,\Phi_{01}+4\,\Phi_{02}
+\Phi_{03}+9\,\Phi_{12}+4\,\Phi_{13}+10\,\Phi_{23}+32}}{\sqrt{70}}\,r 
\leq \Delta \,,
\nonumber
\end{align}
for $m=1,2,3$, respectively, which follow from straightforward computations 
using the matrices in Example~\ref{exm:matrixM} upon setting $\Phi_{ij}:=
\cos(\phi_i-\phi_j)$ for brevity. In each case, the factors multiplying $r$ 
in (\ref{eq-cond-for-m2}) are bounded from above by 1.

\subsection{Preservation of curve end tangent directions}

Although the Bernstein form is more involved in terms of strictly 
satisfying the bound $\|\delta\w\|=\Delta$, it provides a simple means to 
preserve the directions of the curve end derivatives $\r'(0)=\w_0^2$ and 
$\r'(1)=\w_m^2$ by, for example, choosing $\phi_0=\arg(\w_0)$ and $\phi_m
=\arg(\w_m)$, leaving $\phi_1,\ldots,\phi_{m-1}$ and $r_0, \ldots, r_m$ as free parameters 
--- subject to \eqref{ineq} and \eqref{eq:equalPer-B-gen} --- to manipulate the 
curve shape. On the other hand, with the Legendre form and the perturbations 
\eqref{eq:equalPer-L-gen}, the equality $\|\delta\w\|=\Delta$ can be simply satisfied by, e.g., 
choosing $\rho_0=\cdots=\rho_m =\Delta/\sqrt{m+1}$, but the analogous method for
preserving the end derivative directions incurs the complicated conditions
\begin{equation} 
\label{eq-boundary-cond-L}
\arg(\w_0) = \arg\left[\,\sum_{k=0}^m M_{m,0k}\, 
\re^{\ri\,\varphi_k}\right] , \;\;
\arg(\w_m) = \arg\left[\,\sum_{k=0}^m M_{m,mk}\,
\re^{\ri\,\varphi_k}\right] . 
\end{equation}
The following example illustrates these considerations.

\begin{figure}[htb]
\begin{minipage}{0.50\textwidth}
\centering
\includegraphics[width=1\textwidth]{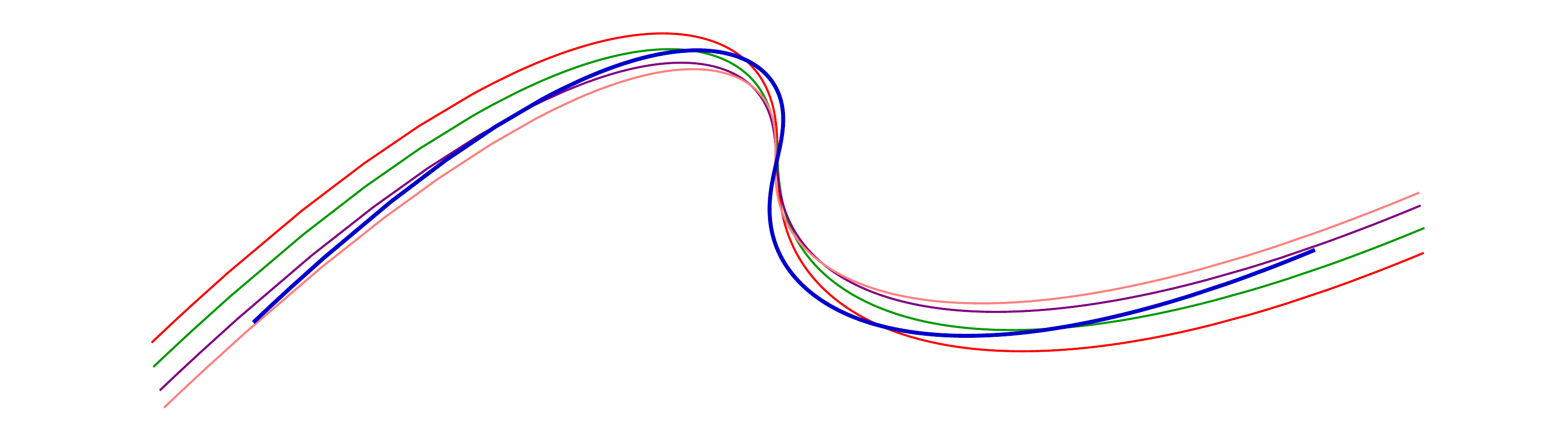}
\end{minipage}
\begin{minipage}{0.50\textwidth}
\centering
\includegraphics[width=1\textwidth]{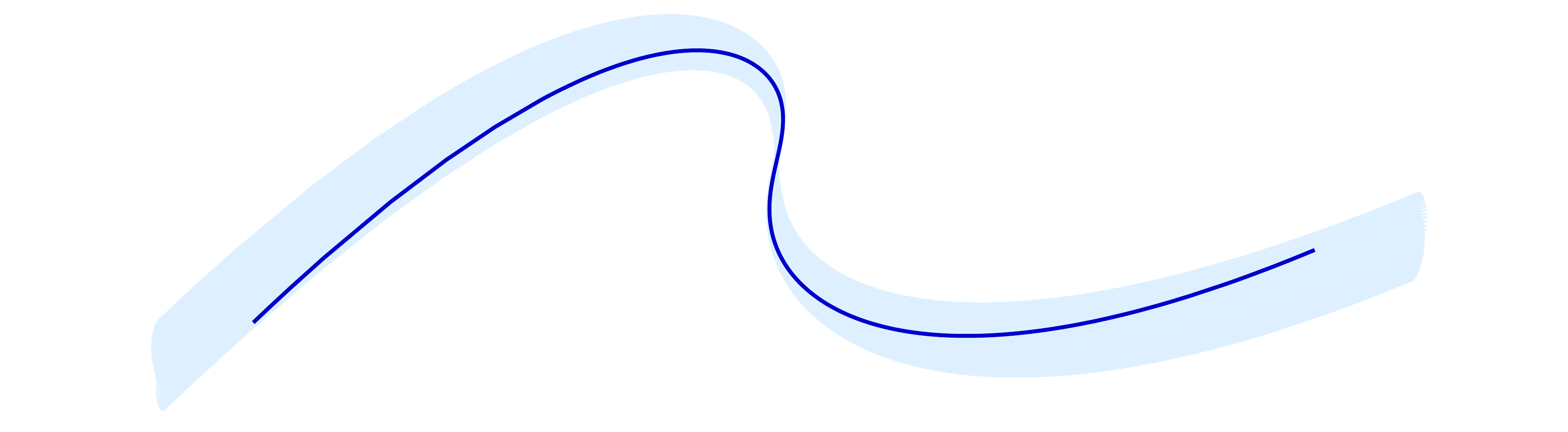}
\end{minipage}
\caption{Left: The prescribed quintic PH curve (blue), with four instances 
modified using the Bernstein basis (different colors), whose pre--images 
satisfy $\|\delta\w\|=0.25$, as described in Example~\ref{exm:example1}. 
Right:  The envelope of the family of all perturbed curves with preserved 
end tangent directions for $r=0.25$.}
\label{fig:example3-1}
\end{figure}

\begin{figure}[htb]
\begin{minipage}{0.50\textwidth}
\centering
\includegraphics[width=1\textwidth]{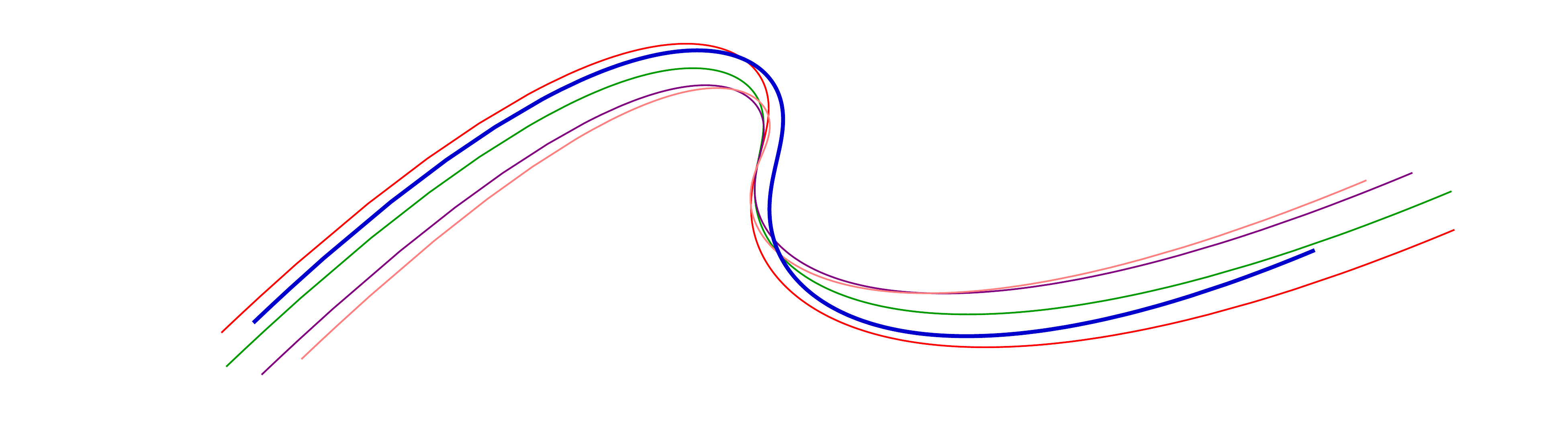}
\end{minipage}
\begin{minipage}{0.50\textwidth}
\centering
\includegraphics[width=1\textwidth]{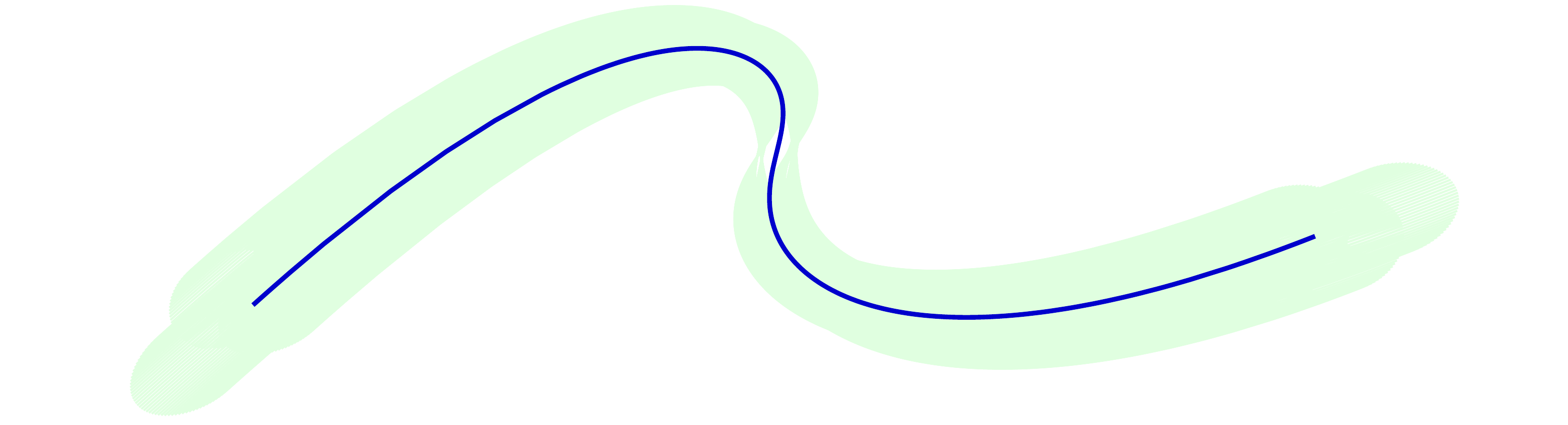}
\end{minipage}
\caption{Left: The prescribed quintic PH curve (blue) with four instances 
modified using the Legendre basis (different colors), whose pre--images 
satisfy $\|\delta\w\|=0.25$, as described in Example~\ref{exm:example1}. 
Right:  The envelope of the family of all perturbed curves with preserved 
end tangent directions for $\rho=0.25/\sqrt{3}$.}
\label{fig:example3-2}
\end{figure}

\begin{exm}
\label{exm:example1}
{\rm 
Consider the quadratic pre--image polynomial $\w(t)$ specified by Bernstein  
coefficients
$$
\w_0 = 5+2\,\ri, \quad \w_1 = -\,3-4\,\ri, \quad \w_2 = 5 + \ri \,,
$$
with corresponding Legendre coefficients
$$
\c_0 = \frac{7}{3}-\frac{1}{3}\,\ri, \quad
\c_1 = -\,\frac{\sqrt{3}}{6}\,\ri, \quad
\c_2 = \frac{8\sqrt{5}}{15}+\frac{11\sqrt{5}}{30}\,\ri \,,
$$
on which we impose perturbations of the form \eqref{eq:equalPer-B-gen} and 
\eqref{eq:equalPer-L-gen} for $m=2$ with equal magnitudes $r=r_0=r_1=r_2$ and $\rho=\rho_0=\rho_1=\rho_2$, satisfying \eqref{ineq} with equality and 
$\Delta=0.25$. With the B\'ezier representation, the end tangent directions are preserved
by choosing $\phi_0=\arg(\w_0)=\arctan(2/5)$ and $\phi_2=\arctan(1/5)$. For 
$\phi_1=0$, $\pi/4$, $\pi/2$, $3\pi/4$ we obtain from \eqref{eq-cond-for-m2} the 
$r$ values for which $\|\delta\w\|=0.25$ as 
$r=0.25245$, $0.25661$, $0.29620$, 
$0.39083$, respectively. Figure~\ref{fig:example3-1} depicts the original 
and four modified PH quintics, all satisfying $\|\delta\w\| =0.25$ --- 
their distances from the original PH quintic are $0.29691,0.30096,0.29884,
0.28626$, respectively. Also shown is the envelope of the family of all 
possible perturbed curves with the prescribed end tangents, for $r=0.25$. 
Note that all the curves have been shifted so that the centroids of their 
B\'ezier control points are coincident. 

With the Legendre representation, the bound $\|\delta\w\|=0.25$ is 
attained, for any angles $\varphi_0$, $\varphi_1$, $\varphi_2$, if and only if 
$\rho=0.25/\sqrt{3}$. To also preserve the end tangent directions, these 
angles must be chosen (see \eqref{eq-boundary-cond-L}) so as to satisfy
\begin{align*}
& 
\arg\,(\re^{\ri\,\varphi_0}-\sqrt{3}\,\re^{\ri\,\varphi_1}
+\sqrt{5}\,\re^{\ri\,\varphi_2}) \,=\, \arctan(2/5) \,, \\
&
\arg\,(\re^{\ri\,\varphi_0}+\sqrt{3}\,\re^{\ri\,\varphi_1}
+\sqrt{5}\,\re^{\ri\,\varphi_2}) \,=\, \arctan(1/5) \,.
\end{align*}
For each of the values $\varphi_0=0,\pi/4,\pi/2,3\pi/4$, four distinct
$(\varphi_1,\varphi_2)$ solutions were identified, defining four different 
perturbations of the quintic PH curve. Figure~\ref{fig:example3-2} compares 
the original PH curve with a representative perturbed PH quintic from each 
of the sets of four solutions. The modified PH quintics have distances 
$0.19350$, $0.22553$, $0.23572$, $0.22451$ from the original PH curve. The 
envelope of the family of modified curves for all $\varphi_0$ values and all solutions is also 
shown (all the curves are shifted so the centroids of their B\'ezier 
control points are coincident).}
\end{exm}

In the preceding discussion, the perturbations incur a \emph{global} change
in the curve. In particular, the curve end points change, which may be 
undesirable in common design contexts. Perturbations to the pre--image 
polynomial that preserve the curve end points are addressed next.

\subsection{Preservation of curve end points}
\label{subsec:pre-image-endpoints}

To eliminate non--essential freedoms, it is customary to consider 
construction of PH curves in \emph{canonical form} \cite{farouki16,farouki22} 
such that $\r(0)=0$ and $\r(1)=1$. The mapping of a PH curve with prescribed 
end points to and from canonical form can be achieved using a simple 
translation/rotation/scaling transformation. Thus, we confine our attention 
to canonical--form PH curves in investigating perturbations that preserve 
the curve end points. Taking $\r(0)=0$ by choice of the integration constant, 
$\r(1)=1$ is achieved through the condition
\be
\label{cf-new}
\int_0^1 \r'(t) \; \rd t \,=\, \int_0^1 \w^2(t) \; \rd t
\,=\, \r(1)-\r(0) \,=\, 1 \,.
\ee
In the Bernstein and Legendre representations \eqref{w} and \eqref{wL-m} 
of the pre--image polynomial $\w(t)$, we set $\boldsymbol{C}=(\c_0,\ldots,
\c_m)^T$ and $\boldsymbol{W}=(\w_0,\ldots,\w_m)^T$. As in \eqref{vecWandC},
the connection between these coefficients is defined by $\boldsymbol{W}= 
M_m\boldsymbol{C}$ (see Lemma~\ref{lemmaCD}). It is easy to see that, for the Legendre form \eqref
{wL-m}, the constraint \eqref{cf-new} reduces to 
\be
\label{eq-point_interpolation-L}
\|\boldsymbol{C}\|_2^2 = \sum_{k=0}^m \c_k^2 \,=\, 1 \,,
\ee
i.e., the Legendre coefficients correspond to points on the unit sphere in 
$\CC^{m+1}$. Setting $\boldsymbol{C}=\boldsymbol{C}_{\!R}+\ri\,\boldsymbol{C}
_{\!I}$, where $\boldsymbol{C}_{\!R},\boldsymbol{C}_{\!I}\in\RR^{m+1}$, we note 
that equation \eqref{eq-point_interpolation-L} is satisfied if and only if
$$ 
\|\boldsymbol{C}_{\!R}\|_2^2-\|\boldsymbol{C}_{\!I}\|_2^2 \,=\, 1 
\quad \mbox{and} \quad \boldsymbol{C}_{\!R}^T\,\boldsymbol{C}_{\!I} \,=\, 0 \,. 
$$
From $\boldsymbol{C}=M_m^{-1}\boldsymbol{W}$, the constraint \eqref{cf-new} 
expressed in terms of the Bernstein coefficients becomes   
\be
\label{eq-point_interpolation-B}
\sum_{k=0}^m \c_k^2 \,=\, \boldsymbol{C}^T\boldsymbol{C} \,=\, 
(M_m^{-1} \boldsymbol{W})^T (M_m^{-1}\boldsymbol{W}) \,=\,  
\boldsymbol{W}^T (M_m^{-1})^T M_m^{-1}\boldsymbol{W} \,=\, 1.
\ee
Setting $G_m=(M_m^{-1})^T M_m^{-1}$, we obtain for $m=1,2,3$ the $(m+1)\times
(m+1)$ matrices with elements $g_{m,jk}$ for $0 \le j,k \le m$ as
\begin{equation*}
G_1 = \frac{1}{6}
\left[
\begin{array}{cc}
 2 & 1 \\
 1 & 2 \\
\end{array}
\right], \quad
G_2=\frac{1}{30} 
\left[
\begin{array}{ccc}
 6 & 3 & 1 \\
 3 & 4 & 3 \\
 1 & 3 & 6 \\
\end{array}
\right], \quad     
G_3=\frac{1}{140} 
\left[
\begin{array}{cccc}
 20 & 10 & 4 & 1 \\
 10 & 12 & 9 & 4 \\
 4 & 9 & 12 & 10 \\
 1 & 4 & 10 & 20 \\
\end{array}
\right],
\end{equation*}
and in the cases $m=1,2,3$ equation \eqref{eq-point_interpolation-B} then 
reduces to 
\begin{gather*}
 \w_0^2+\w_1^2+\w_1 \w_0  \,=\, 3, \\
 3\,\w_0^2+3\,\w_2^2+2\,\w_1^2+3\,(\w_0+\w_2)\,\w_1+\w_0\w_2 \,=\, 15 \,,  \\
 10\,(\w_0^2+\w_3^2)+6\,(\w_1^2+\w_2^2)+10\,(\w_0\w_1+\w_2\w_3) 
+\; 4\,(\w_2\w_0+\w_1\w_3)+\w_3\w_0+9\,\w_1\w_2  \,=\, 70 \,.
\end{gather*}

To ensure that the conditions \eqref{eq-point_interpolation-L} and \eqref
{eq-point_interpolation-B} are fulfilled upon substituting $\c_k\to\c_k 
+\delta\c_k$ and $\w_k\to\w_k+\delta\w_k$ for $k=0,\dots, m$, the coefficients 
$\delta\c_k$ and $\delta\w_k$ must satisfy  
\be
\label{endCondL}
\sum_{k=0}^m \delta\c_k^2+2\sum_{k=0}^m \c_k\,\delta\c_k \,=\,
\delta\boldsymbol{C}^T(\delta\boldsymbol{C}+2\,\boldsymbol{C}) 
\,=\, 0 \,,
\ee
and 
\be
\label{endCondB}
\delta\boldsymbol{W}^TG_m (\delta\boldsymbol{W}+2\,\boldsymbol{W}) \,=\,
\sum_{j=0}^m \sum_{k=0}^m g_{m,jk}\,\delta\w_j\,(\delta\w_k+2\,\w_k) \,=\, 0 \,,
\ee
for the prescribed $\c_k$ and $\w_k$ values satisfying \eqref
{eq-point_interpolation-L} and \eqref{eq-point_interpolation-B}. In addition,
the perturbations must satisfy the bounds $\|\delta\boldsymbol{C}\|_2\leq
\Delta$ and $\|M_m^{-1}\delta\boldsymbol{W}\|_2\leq\Delta$. 
Writing $\delta\boldsymbol{C}=\delta\boldsymbol{C}_{\!R}+\ri\,\delta
\boldsymbol{C}_{\!I}$, the condition (\ref{endCondL}) is equivalent to two 
scalar equations in the real vectors $\delta\boldsymbol{C}_{\!R}$ and 
$\delta\boldsymbol{C}_{\!I}$, namely
\begin{equation} \label{endCondL-1}
\begin{split}
\|\delta\boldsymbol{C}_{\!R}\|_2^2-\|\delta\boldsymbol{C}_{\!I}\|_2^2+
2\,\delta\boldsymbol{C}_{\!R}^T\boldsymbol{C}_{\!R}-
2\,\delta\boldsymbol{C}_{\!I}^T\boldsymbol{C}_{\!I} & \,=\, 0 \,,\\
\delta\boldsymbol{C}_{\!R}^T\delta\boldsymbol{C}_{\!I}+
\delta\boldsymbol{C}_{\!R}^T\boldsymbol{C}_{\!I}+
\delta\boldsymbol{C}_{\!I}^T\boldsymbol{C}_{\!R} & \,=\, 0 \,.
\end{split}
\end{equation}
Writing $\delta\boldsymbol{W}=\delta\boldsymbol{W}_{\!R}+\ri\,\delta\boldsymbol
{W}_{\!I}$, we have $\delta\boldsymbol{W}_{\!R}=M_m\,\delta\boldsymbol{C}_{\!R}$ 
and $\delta\boldsymbol{W}_{\!I}=M_m\,\delta\boldsymbol{C}_{\!I}$, so these 
equations can also be expressed in terms of $\delta\boldsymbol{W}_{\!R}$ and 
$\delta\boldsymbol{W}_{\!I}$.

\begin{figure}[htbp]
\centerline{
\includegraphics[width=0.4\textwidth]{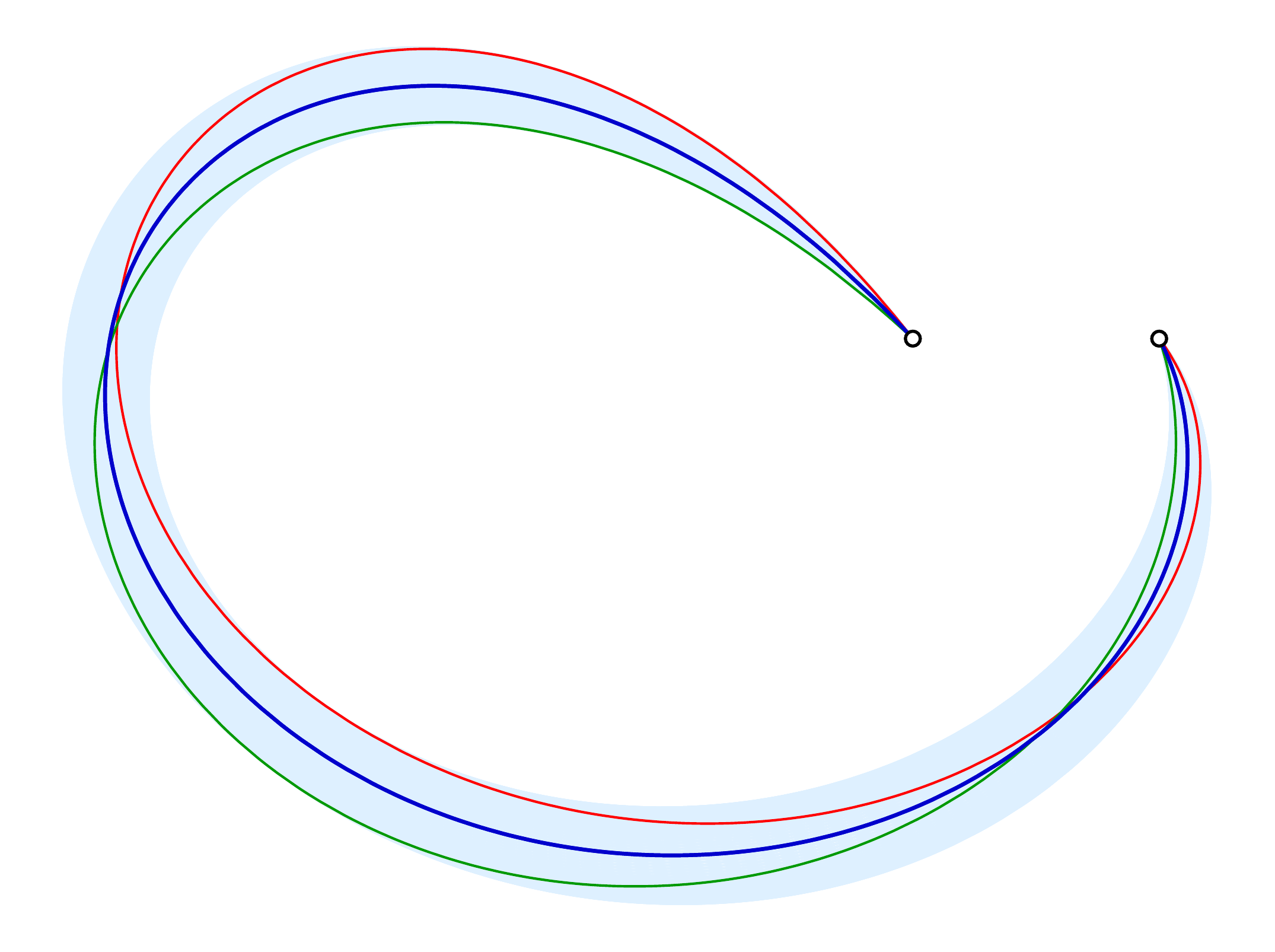} \qquad
\includegraphics[width=0.4\textwidth]{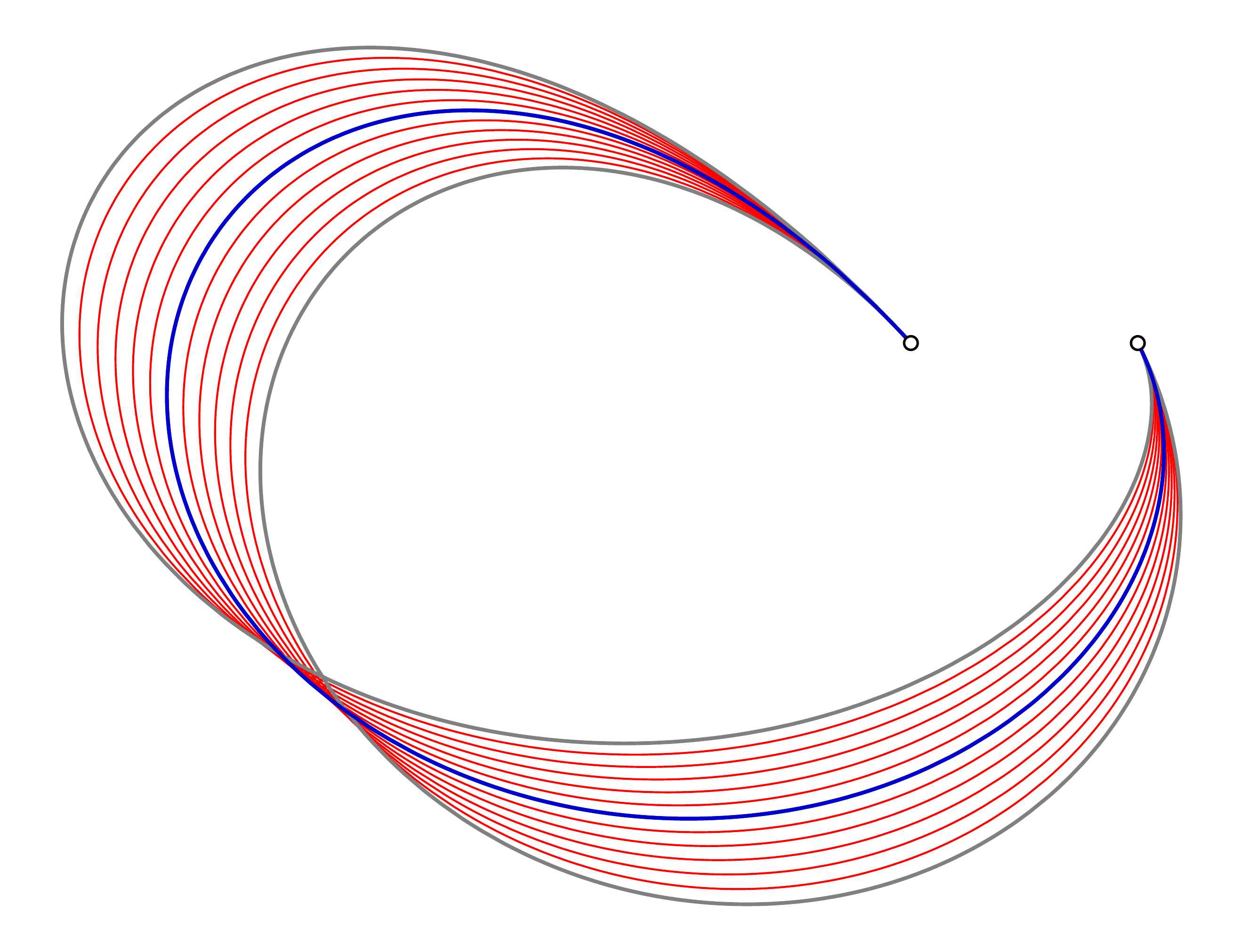}
}
\caption{Left: The perturbed PH quintics (red and green) for $\|\delta \w\| = 0.1$, with the same end points as the quintic PH curve (blue) in Example~\ref{exm:example5}, 
defined by $\varphi=0$. The envelope of all the solutions for $\varphi\in (-\,\pi,\pi\,]$ is also shown (light blue).
Right: The modified PH quintic curves that preserve end points and end tangent directions, as described in Example~\ref{exm:tangents}, as $\|\delta \w\|$ varies. 
} 
\label{fig:example5}
\end{figure}

In general, the identification of perturbations $\delta\w$ that maintain
the perturbed curve in canonical form and satisfy $\|\delta\w\|=d\leq
\Delta$ for some chosen $d$ entails the solution of a rather complicated
non--linear system. We now propose a simple sufficient way to determine such
perturbations, focusing on the Legendre representation with coefficients
$\delta \c_k$ of the form \eqref{eq:equalPer-L-gen}. The idea is to consider
coefficients for which the angles are the same, i.e., $\varphi:=\varphi_0=
\ldots = \varphi_m$. Then 
\begin{align*}
  \|\delta \VCR\|_2^2 - \|\delta \VCI\|_2^2
  \,=\, \sum_{k=0}^m \rho_k^2 \cos{(2\varphi)} \,, \quad
  \delta \VCR^T \delta \VCI \,=\, \frac{1}{2}
  \sum_{k=0}^m \rho_k^2 \sin{(2\varphi)} \,,
\end{align*}
so the equations \eqref{endCondL-1}, together with the condition 
\begin{equation} \label{eq:cond-norm-w-d}
\|\delta \w \|^2 \,=\, \sum_{k=0}^m \rho_k^2 \,=\, d^2 \,,
\end{equation}
simplify to the two linear equations
\begin{equation} \label{eq-perturbation-2}
\begin{split}
d^2 \cos{(2\varphi)} + 2 \sum_{k=0}^m \rho_k (\cos(\varphi) \, \mbox{Re}(\c_k) - 
\sin(\varphi) \, \mbox{Im}(\c_k)) & \,=\, 0 \,,\\  
d^2 \sin{(2\varphi)} + 2 \sum_{k=0}^m \rho_k (\cos(\varphi) \, \mbox{Im}(\c_k) + 
\sin(\varphi) \, \mbox{Re}(\c_k)) & \,=\, 0 \,,
\end{split}
\end{equation}
for the radii $\rho_k$, $k=0,1, \dots, m$, that depend on the chosen angle
$\varphi$. 
Thus, the only non--linear aspect is to satisfy \eqref{eq:cond-norm-w-d}.
The following example demonstrates this construction.
\begin{exm}
\label{exm:example5}
{\rm Consider the quadratic pre--image polynomial with Legendre coefficients 
$$
\c_0 = 2 -  \ri, \quad
\c_1 = 1 + 2 \ \ri, \quad
\c_2 = -1 + 0 \ \ri,
$$
that satisfy \eqref{eq-point_interpolation-L}. We choose $d=0.1$ for
perturbations of the form 
\eqref{eq:equalPer-L-gen} with $\varphi_0=\cdots = \varphi_m = \varphi$ for any $\varphi \in (-\,\pi,\pi\,]$. 
From \eqref{eq-perturbation-2} we obtain
$$
\rho_1= \frac{\rho_0}{2}-\frac{\sin (\varphi)}{400}, \quad
\rho_2 = \frac{5\rho_0}{2}-\frac{\sin (\varphi)}{400}+\frac{\cos(\varphi)}{200},
$$
and \eqref{eq:cond-norm-w-d} reduces to a quadratic equation for $\rho_0$,
namely
$$
\frac{15}{2}\, \rho_0^2 + \frac{5 \cos (\varphi )\, -3 \sin (\varphi )}{200} \, \rho_0 +\frac{\cos (2 \varphi)-2 \sin (2 \varphi )-1597}{160000} = 0,
$$
which has a positive discriminant 
for any $\varphi$. Thus, there are always two admissible perturbations $\delta \w$, shown in Figure~\ref{fig:example5} (left) for the choice $\varphi=0$  (red and green curves), together with the original (blue) curve, and the envelope of all solutions for
$\varphi\in (-\,\pi,\pi\,]$.
}
\end{exm}

With the Legendre representation it is easy to construct perturbations $\delta \w(t)$ that preserve end points of the given curve, while the B\'ezier representation is more convenient for preserving end tangent directions. The constraints on the coefficients of $\delta\w(t)$ imposed by preserving end points and end tangent directions, coupled with the non--linear dependence of $\|\delta\w\|$ on those coefficients, makes it difficult to formulate schemes that guarantee an \emph{a priori} satisfaction of the bound \eqref{ineq}. As 
a practical solution for PH quintics, we consider here coefficients $\delta \w_0=r\exp(\ri\,\phi_0)$ and $\delta\w_2=r\exp(\ri\,\phi_2)$, where $\phi_0=\arg(\w_0)$ and $\phi_2=\arg(\w_2)$, for a prescribed $r$ value, to preserve 
the end tangent directions. Preservation of the end points can then be achieved by solving the $m=2$ instance of equation \eqref{endCondB}
for $\w_0,\w_1,\w_2$ values that define a canonical--form PH quintic, as a 
quadratic equation in $\delta\w_1$. Since this incurs modest computational 
effort, it is amenable to real--time user modification of $r$ to ensure 
satisfaction of \eqref{ineq}. However, adding the constraint $\|\delta \w\| = d$ for some $d\leq \Delta$ adds 
one non--linear equation, but since the whole non--linear system is algebraic, it is possible to compute all the solutions using a computer algebra system.

\begin{figure}[htbp]
\begin{minipage}{0.50\textwidth}
\centering
\includegraphics[width=1\textwidth]{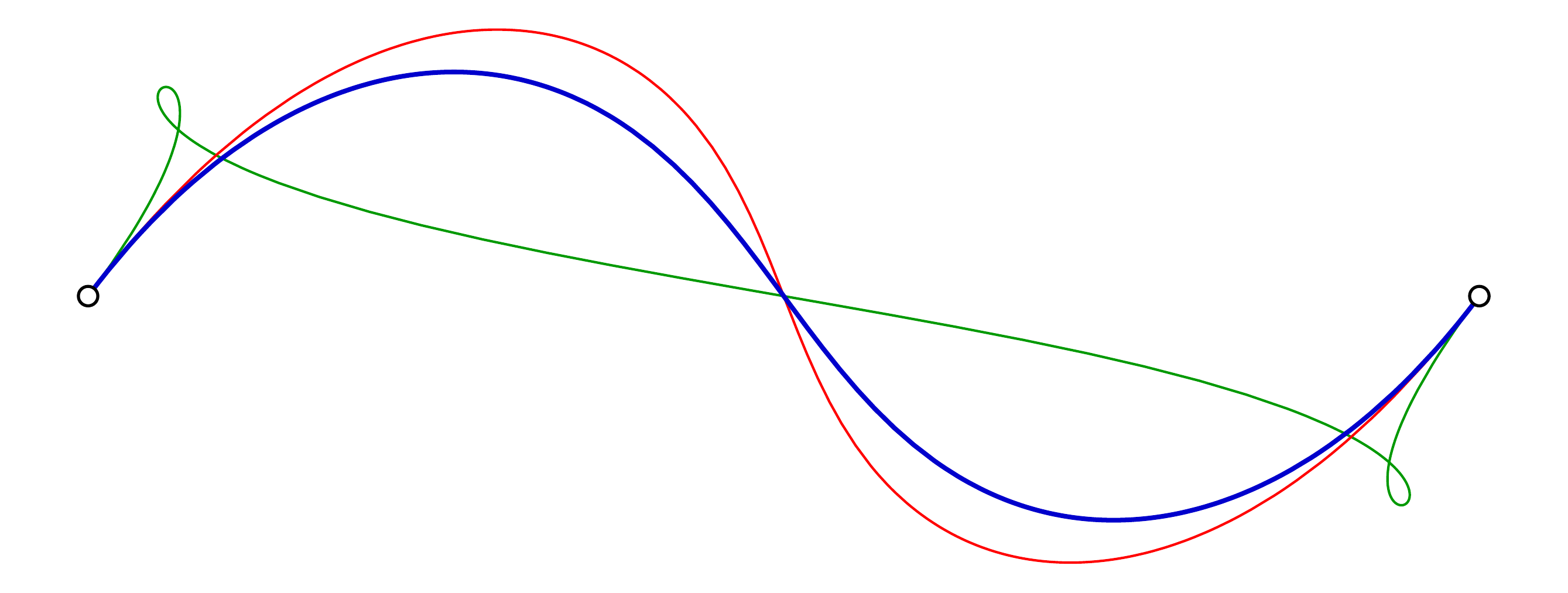}
\end{minipage}
\begin{minipage}{0.50\textwidth}
\centering
\includegraphics[width=1\textwidth]{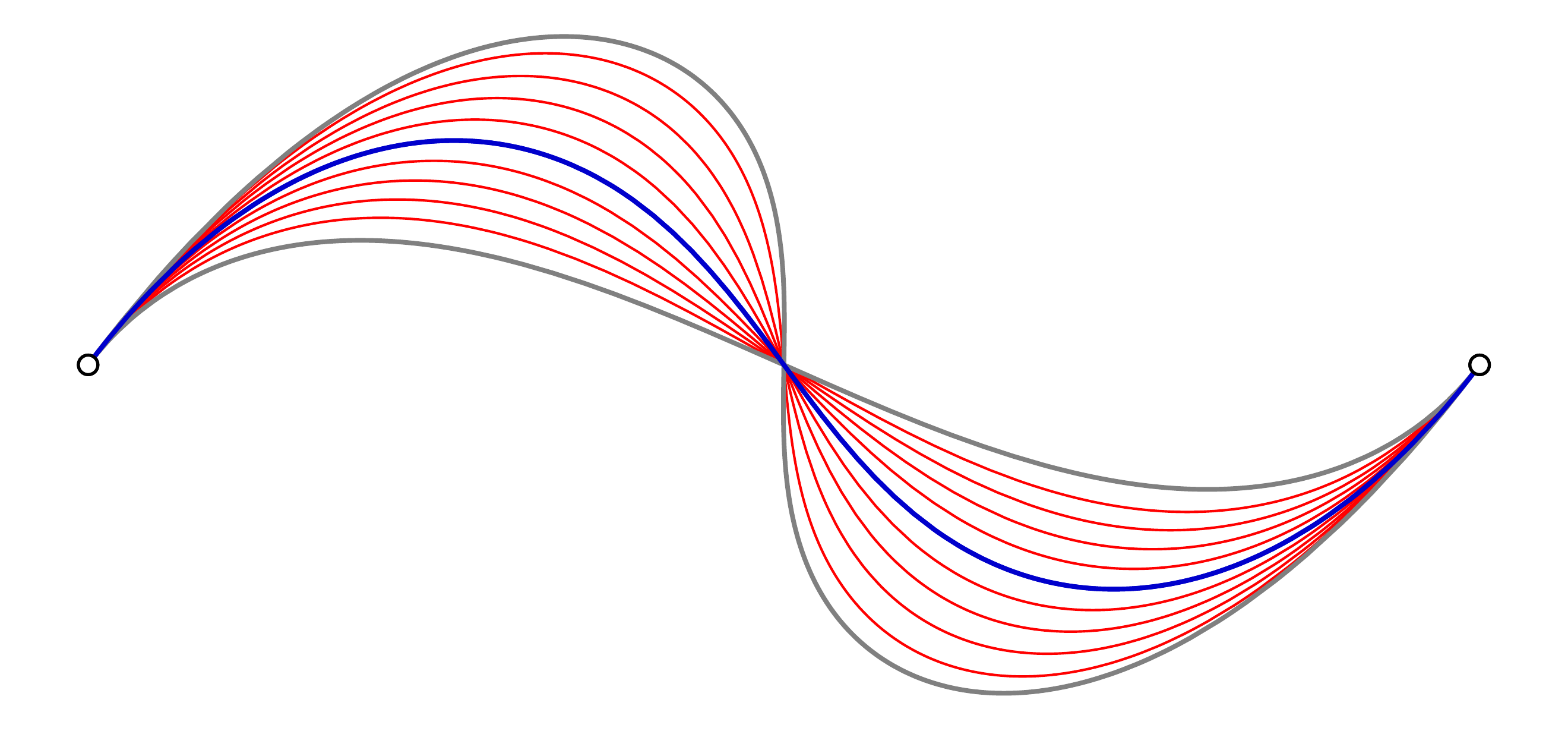}
\end{minipage}
\caption{Left: The given canonical--form PH quintic (blue) with two modifications 
(red, green) that preserve the end points and tangent directions for a fixed $r=0.2$, as described in 
Example~\ref{exm:tangents}.
Right: The modified curves satisfying $\|\delta \w\| \leq 0.25$ as $r$ varies. 
}
\label{fig:example6}
\end{figure}

\begin{exm}
\label{exm:tangents}
{\rm Consider a canonical--form PH quintic defined by a quadratic 
pre--image polynomial with Bernstein coefficients
\begin{align*}
\w_0 &= \w_2 = \sqrt{2} + \frac{\sqrt{2}}{2} \,\ri  
\,, \\
\w_1 & =  \frac{\sqrt{5 (9+\sqrt{97})}-6\sqrt{2}}{4}-
\frac{\sqrt{-27+5 \sqrt{97}+6 \sqrt{10 (\sqrt{97}-9)}}}{4}\,\ri
\,.
\end{align*}
To preserve the end tangent angles $\theta_0=\arg(\w_0)$ and $\theta_2=\arg(\w_2)$, we set 
$\delta\w_0=r\exp(\ri \,\theta_0)$, $\delta\w_2= r\exp(\ri \,\theta_2)$ for some chosen $r$, and compute $\delta \w_1$
from equation \eqref{endCondB} for $m = 2$. With $r=0.2$ we obtain two solutions 
\[
\delta\w_1 \,=\, -\,0.33476348-0.29109547\,\ri \,, 
\quad \delta\w_1 \,=\, -\,5.05586773+1.05285093\,\ri \,.
\]
The corresponding $\|\delta\w\|$ values are $0.102659$, $1.802944$.
Figure~\ref{fig:example6} (left) shows the resulting curves. The first solution (red curve) is 
evidently a very reasonable modification of the original curve (blue), preserving 
its end points and end tangents. Although the second solution (green curve) also has 
this property, it exhibits tight loops --- a common feature \cite
{farouki16,farouki95} among the multiple solutions to PH quintics that
satisfy given constraints --- and is discarded on the basis of the large
$\|\delta\w\|$ value. The perturbed curves with $\|\delta \w\| \leq 0.25$ for 
$r = -0.4, -0.3, \dots, 0.3, 0.4$ are shown (red curves) in Figure~\ref{fig:example6} (right)
together with the two curves (gray) having $\|\delta \w\| =0.25$, obtained for 
$r= -0.52962446$ and $r=0.47220859$.

Choosing the pre-image polynomial defined in Example~\ref{exm:example5}, with  
Bernstein coefficients 
$$
\w_0 = (2-\sqrt{3}-\sqrt{5}) -(1+2 \sqrt{3}) \,\ri\,, \; 
\w_1 = 2 \left(1+\sqrt{5}\right) - \ri\,, \;
\w_2 = (2+\sqrt{3}-\sqrt{5}) + (2 \sqrt{3}-1)\, \ri
$$
and the choices $r=-0.5,-0.4,\ldots,0.4,0.5$, the modified quintic PH curves that preserve end points and end tangent directions, and satisfy $\|\delta \w\| \leq 0.25$  
are shown (red curves) in Figure~\ref{fig:example5} (right) together with the two curves (gray) having $\|\delta \w\| =0.25$, obtained for 
$r= -0.59313245$ and $r=0.60204179$.
}
\end{exm}

\section{Modification of PH curve arc lengths}
\label{sec:arclen}

The total arc length $S$ of a planar PH curve $\r(t)$ is intimately related 
to the norm of its pre--image polynomial $\w(t)$, since
\be
S \,=\, \int_0^1 \sigma(t) \, \rd t \,=\, 
\int_0^1 |\w(t)|^2 \, \rd t \,=\, \|\w\|^2 \,.\nonumber
\ee
The arc length $S$ can be changed by a specified amount $\delta S>0$ 
by choosing $\delta\w(t)$ to have the norm $\sqrt{\delta S}$ and to be 
orthogonal to $\w(t)$, since
\[
\int_0^1 |\w(t)+\delta\w(t)|^2 \, \rd t \,=\,
\|\w\|^2+\|\delta\w\|^2+2\,\mbox{Re}(\langle\w,\delta\w\rangle) \,,
\]
where $\|\w\|^2=S$, $\|\delta\w\|^2=\delta S$, and $\mbox{Re}(\langle\w,
\delta\w\rangle)=0$ when $\w(t)$ and $\delta\w(t)$ are orthogonal. 
We focus here on the Legendre representation 
\[
\w(t) \,=\, \sum_{k=0}^m \c_kL_k(t)
\quad \mbox{and} \quad
\delta\w(t) \,=\, \sum_{k=0}^m \delta\c_kL_k(t) 
\]
of the pre--image polynomial and its perturbation. Recalling the notation from
Section~\ref{subsec:pre-image-endpoints} we set $\c_k=c_{k,1}+\ri\, c_{k, 2}$
and form, as in Section~\ref{sec:ocurves}, the real vectors 
\begin{align}  \label{eq:aAndq}
\a  \,=\, (c_{0,1},c_{0,2}, c_{1,1}, c_{1,2}, \dots, c_{m,1}, c_{m,2})^T  
\quad \mbox{and} \quad \g  \,=\, \a+\mbox{sign}(c_{0,1})\|\a\|_2 (1,0, \dots, 0)^T \,.
\end{align}
The second through last columns 
of the $(2m+2)\times(2m+2)$ matrix $Q=(q_{j,k})_{j,k=1}^{2m+2}$ defined in terms of 
$\a$ and $\g$ defined in \eqref{eq:aAndq} by the formula \eqref{Q} then identify the coefficients of the polynomials
\be
\b_{k}(t) \,=\, \sum_{j=0}^{m} \b_{k,j} L_{j}(t)\,, \quad \b_{k, j}:= q_{2j+1,k+1} + \ri\, q_{2j+2,k+1},  
\quad k=1,2,\dots, 2m+1,\nonumber
\ee
that form the orthonormal basis for degree $m$ complex polynomials orthogonal to $\w(t)$ with norms 
$\|\b_k\|=1$, $k=1,2, \dots, 2m+1$. Thus, for any real values $\gamma_1, \gamma_2, \dots, \gamma_{2m+1}$ a 
perturbation polynomial of the form
\be
\label{dwgamma}
\delta\w(t) \,=\, \sum_{k=1}^{2m+1} \gamma_k \b_{k}(t)
\ee
is orthogonal to $\w(t)$ and has norm $\|\delta\w\|= \sqrt{\gamma_1^2+\gamma_2^2+\dots + \gamma_{2m+1}^2}$. Thus, by assigning values $\gamma_1,\gamma_2, \dots, \gamma_{2m+1}$ 
that satisfy 
\be
\label{dSeqn}
\gamma_1^2+\gamma_2^2+\dots + \gamma_{2m+1}^2 \,=\,\delta S \,,
\ee
the perturbed pre--image polynomial $\w(t)+\delta\w(t)$ generates a PH curve
with arc length $S+\delta S$. 

To ensure that the modified curve has the same end points as the given PH
curve $\r(t)$, assumed to be in canonical form, $\delta\w(t)$ must also
satisfy the condition \eqref{endCondL}, which can be reduced to the quadratic equation
\begin{equation}
\label{Qeqn}
\sum_{j,k = 1}^{2m+1} \f_{j,k} \gamma_j \gamma_k + \sum_{k=1}^{2m+1} \f_{k} \gamma_k \,=\, 0,
\quad {\rm where} \quad \f_{j,k} := \sum_{\ell=0}^m \b_{j,\ell} \b_{k,\ell}, \quad 
\f_{k} := 2 \sum_{\ell=0}^m \b_{k,\ell} \c_{\ell},
\end{equation}
in $\gamma_1, \dots,\gamma_{2m+1}$. 
Equation \eqref{dSeqn} and the real and imaginary parts of equation \eqref
{Qeqn} constitute a system of three quadratic equations for $2m+1$ factors
$\gamma_1,\dots,\gamma_{2m+1}$ in \eqref{dwgamma}, which allows one to fix
$2m -2$ of them, and then
solve the system using Newton--Raphson iterations or some other algebraic
solver.

\begin{rmk} \label{rem-matrixQ}
{\rm If we denote by $Q_R$ the sub--matrix of $Q$ with rows $1,3, \dots, 2m+1$ and columns 
$2,3,\ldots,2m+2$, and by 
$Q_I$ the sub-matrix of $Q$ with rows $2,4,\dots, 2m+2$ and columns 
$2,3,\ldots,2m+2$, and we
define the complex matrix $\boldsymbol{Q} = Q_R + \ri \, Q_I$, 
the vector $\delta\boldsymbol{C}$ of the coefficients of $\delta \w$ can be expressed as
$\delta\boldsymbol{C} =  \boldsymbol{Q}\,\boldsymbol{\gamma}$ for 
$\boldsymbol{\gamma}=(\gamma_1, \gamma_2, \dots, \gamma_{2m+1})^T$, which gives a more compact representation of \eqref{Qeqn}, namely
$$
\boldsymbol{\gamma}^T \boldsymbol{F} \boldsymbol{\gamma} + 2\,\boldsymbol{\gamma}^T \boldsymbol{Q}\,\boldsymbol{C} = 0 \,, \quad \boldsymbol{F} := \boldsymbol{Q}^T \boldsymbol{Q} \,.
$$
}
\end{rmk}

\begin{figure}[htbp]
\begin{minipage}{0.50\textwidth}
\centering
\includegraphics[width=1\textwidth]{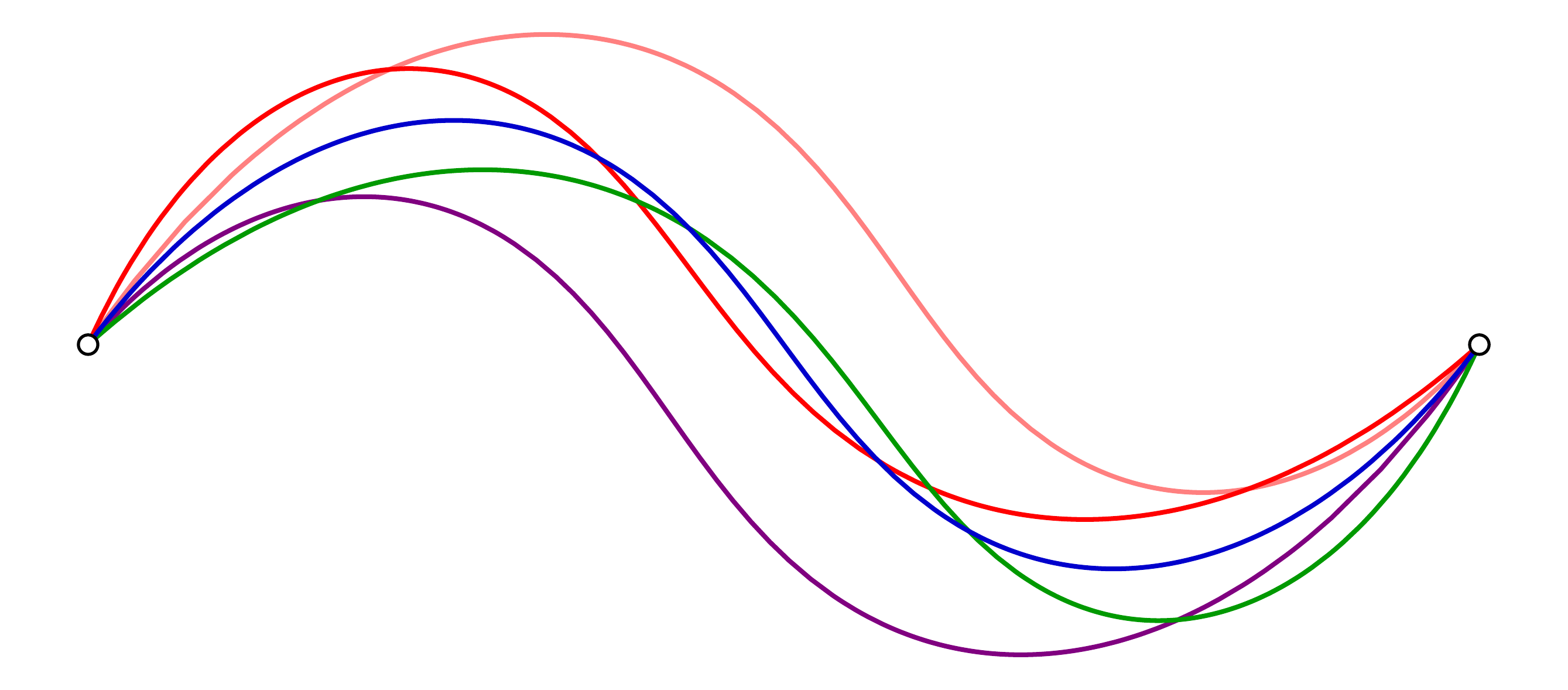}
\end{minipage}
\begin{minipage}{0.50\textwidth}
\centering
\includegraphics[width=1\textwidth]{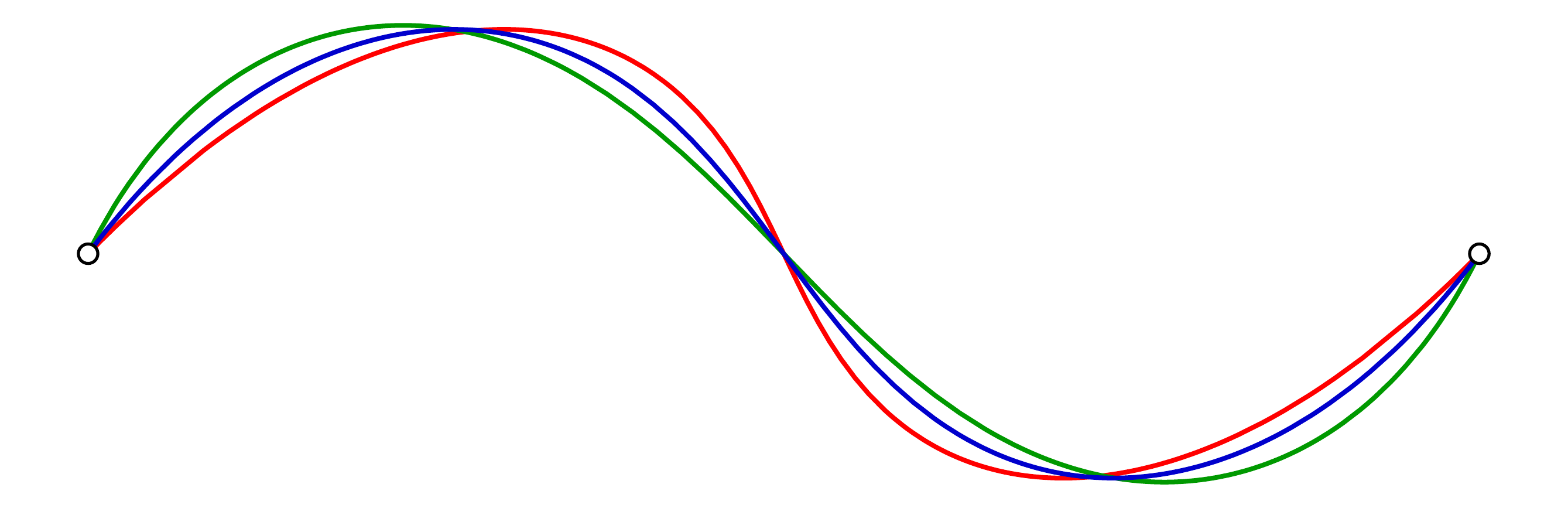}
\end{minipage}
\caption{The PH quintic from Example~\ref{exm:tangents} (blue), together with 
PH quintics (different colors) with arc lengths increased by $\delta S=0.01$, sharing the 
same end points, computed by fixing $\gamma_4=\gamma_5=0$ (left) and 
$\gamma_2=\gamma_3=0$ (right).}
\label{fig:curve12}
\end{figure}

\begin{exm}
{\rm Consider the PH quintic specified by the pre--image polynomial in 
Example~\ref{exm:tangents}, with arc length $S=1.23740482$.
The complex matrix $\boldsymbol{Q}$, defined in Remark~\ref{rem-matrixQ}, is
$$
\boldsymbol{Q}= \left[
\begin{array}{ccccc}
 0.048391+0.998792 \,\ri & 0 & 0 & -0.148557+0.003707 \,\ri & -0.305920+0.007634 \,\ri \\
 0 & 1 & \ri & 0 & 0 \\
 0.003707+0.007634 \,\ri & 0 & 0 & 0.988619-0.023436 \,\ri & -0.023436+0.951738 \,\ri \\
\end{array}
\right] \,.
$$
The perturbation $\delta \w(t)$ in \eqref{dwgamma} is expressed in terms
of the five parameters $\gamma_1\dots, \gamma_5$. Fixing two of them and
choosing $\delta S=0.01$, the remaining parameters may be computed as the
solution of equations \eqref{dSeqn}--\eqref{Qeqn}. 
For $\gamma_4=\gamma_5=0$ four solutions are identified:
\begin{align*}
(\gamma_1,\gamma_2,\gamma_3) & =
(0.0047585271, \pm 0.074073623, \mp 0.067010856)\,,
\\
(\gamma_1,\gamma_2,\gamma_3) & =
(-0.0047585271, \pm 0.074073623, \pm 0.067010856) \,.
\end{align*}
Fixing $\gamma_2=\gamma_3=0$ the sytem has only two solutions: 
\begin{align*}
(\gamma_1,\gamma_4,\gamma_5) & =
(-0.032364637, 0.094555975, 0.0034202026)\,,
\\
(\gamma_1,\gamma_4,\gamma_5) & =
(0.030467676, -0.094859055, 0.0085720767) \,.
\end{align*}
Using these values, the Legendre coefficients of $\delta\w(t)$ follow from $\delta\boldsymbol{C} =  \boldsymbol{Q} \boldsymbol{\gamma}$. 
The resulting PH quintics with increased arc length, generated by the modified 
pre--image polynomials $\w(t)+\delta\w(t)$, are shown in Figure~\ref{fig:curve12}. 
}
\end{exm}

\begin{figure}[htbp]
\centerline{
\includegraphics[width=0.4\textwidth]{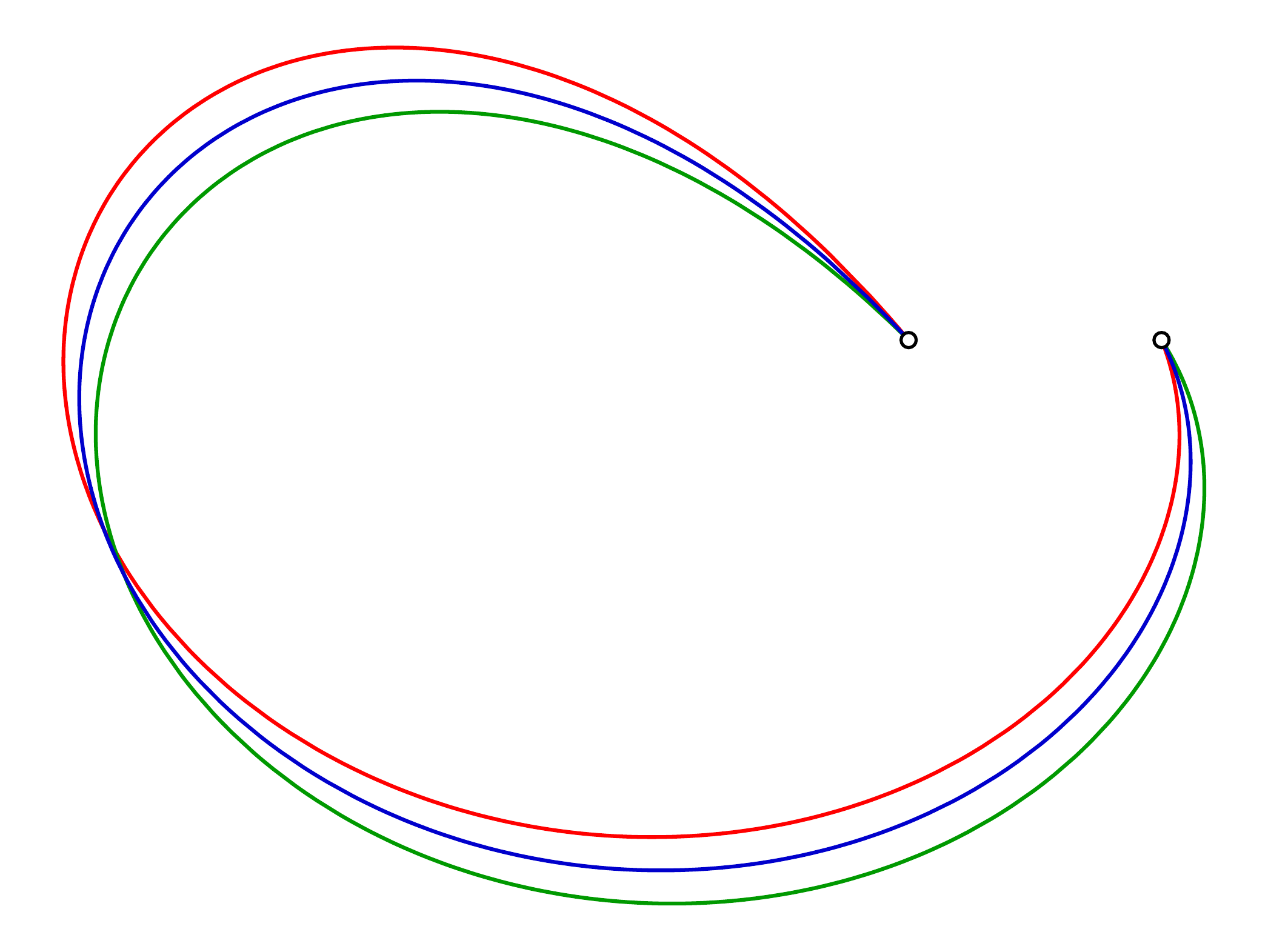} \qquad
\includegraphics[width=0.4\textwidth]{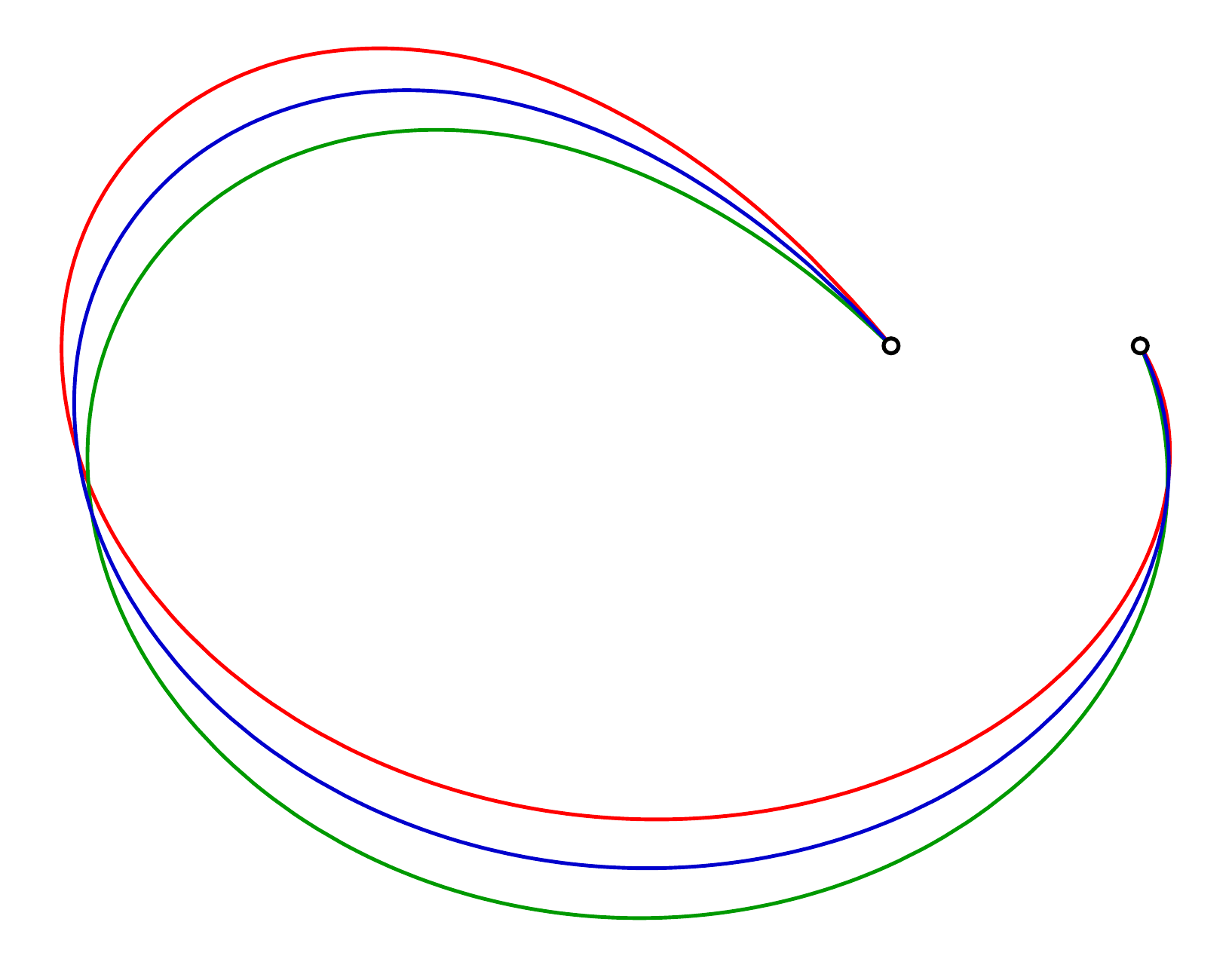}
}
\caption{The PH quintic from Example~\ref{exm:arclen-2} (blue), together with 
PH quintics (red, green) with arc lengths increased by $\delta S=0.01$, sharing the 
same end points, computed by fixing $\gamma_4=\gamma_5=0$ (left) and 
$\gamma_2=\gamma_3=0$ (right).}
\label{fig:exm-8}
\end{figure}

\begin{exm}
\label{exm:arclen-2}
{\rm As a final example we choose the PH quintic specified by the pre--image polynomial in 
Example~\ref{exm:example5}, with arc length $S=11$, and follow the same steps as in the previous example. For the choice $\delta S=0.01$ and $\gamma_4=\gamma_5=0$ there are two solutions,
\begin{align*}
(\gamma_1,\gamma_2,\gamma_3) & =
(-0.068622784, 0.069792544, -0.020491812)\,,
\\
(\gamma_1,\gamma_2,\gamma_3) & =
(0.068681604, -0.069717905, 0.020548745) \,,
\end{align*}
shown in Figure~\ref{fig:exm-8} (left). Fixing $\gamma_2=\gamma_3=0$, we obtain 
\begin{align*}
(\gamma_1,\gamma_4,\gamma_5) & =
(-0.034621641, -0.083559293, -0.042651924)\,,
\\
(\gamma_1,\gamma_4,\gamma_5) & =
(0.031439019, 0.083161950, 0.045778578) \,.
\end{align*}
The corresponding curves are shown in Figure~\ref{fig:exm-8} (right).
}
\end{exm}

\section{Closure}
\label{sec:close}

Interpreting planar polynomial curves as complex--valued functions of a real
parameter $t\in [\,0,1\,]$ facilitates the introduction of an inner product,
norm, and metric function, that permit measurement of curve magnitudes and 
of the distances and angles between curves. The concept of orthogonal curves
is then possible, leading to a procedure to construct a basis spanning all 
planar curves that are orthogonal to a given planar curve.

These concepts were applied to the complex pre--image polynomials that 
define planar Pytha\-gorean--hodograph (PH) curves, to develop schemes that 
allow bounded modifications of a given PH curve, without compromising its PH 
nature. Specializations of these schemes that accommodate preservation of 
curve end points and end tangents have also been presented, and the use of 
an orthogonal basis for a given PH curve pre--image polynomial to achieve 
a desired change in the arc length of the PH curve was demonstrated.

The methodology presented herein may also be generalized to the spatial
Pythagorean--hodograph curves through the quaternion representation \cite
{choi02,farouki02} and preliminary results have already been reported 
in \cite{farouki24}. A further domain of interest concerns a possible 
adaptation of the methodology from curves to parametric surfaces, defined
as vector quaternion polynomial functions of two parameters over triangular 
or rectangular domains.
\vspace{0.5cm}

{\noindent \sl {\bf Acknowledgments}.}
The second and the fourth
author have been partly supported by the research program P1-0288 and the research grants N1-0137, N1-0237 and J1-3005 by the Slovenian Research and Innovation Agency.
The third author has been partly supported by the research program P1-0404 and the research grants N1-0296, N1-0210 and J1-4414 
by the Slovenian Research and Innovation Agency.


\def\AMM{{\it Amer.\ Math.\ Monthly\ }}
\def\ACM{{\it Adv.\ Comp.\ Math.\ }}
\def\ACMTMS{{\it ACM Trans.\ Math.\ Software\ }}
\def\ACMTOG{{\it ACM Trans.\ Graphics\ }}
\def\AMC{{\it Appl.\ Math.\ Comput.\ }}
\def\CAD{{\it Comput.\ Aided Design }}
\def\CAEJ{{\it Comput.\ Aided Eng. J.\ }}
\def\CAGD{{\it Comput.\ Aided Geom.\ Design }}
\def\CAVW{{\it Comput.\ Anim.\ Virt.\ Worlds }}
\def\CG{{\it Computers \& Graphics }}
\def\CVGIP{{\it Comput.\ Vision, Graphics, Image\ Proc.\ }}
\def\GM{{\it Graph.\ Models\ }}
\def\IBMJRD{{\it IBM J.\ Res.\ Develop.\ }}
\def\JCAM{{\it J.\ Comput.\ Appl.\ Math.\ }}
\def\JGG{{\it J.\ Geom.\ Graphics }}
\def\JMAA{{\it J.\ Math.\ Anal.\ Appl.\ }}
\def\JSC{{\it J.\ Symb.\ Comput.\ }}
\def\MC{{\it Math.\ Comp.\ }}
\def\MMA{{\it Math.\ Model.\ Anal.\ }}
\def\MMAS{{\it Math.\ Methods\ Appl.\ Sci.\ }}
\def\NA{{\it Numer.\ Algor.\ }}
\def\NMTMA{{\it Numer.\ Math.\ Theor.\ Meth.\ Appl.\ }}
\def\PAMS{{\it Proc.\ Amer.\ Math.\ Soc.\ }}
\def\SIAMJNA{{\it SIAM J.\ Numer.\ Anal.\ }}
\def\SIAMR{{\it SIAM Rev.\ }}
\def\TASCE{{\it Trans.\ Amer.\ Soc.\ Civil Eng.\ }} 
\def\TF{{\it Transport.\ Forum }}

\begin{flushleft}

\end{flushleft}

\end{document}